\documentclass{article}

\usepackage{amsmath,amsfonts,amsthm}
\usepackage{graphicx}

\newcommand{\nz}{\mathbb{N}}
\newcommand{\gz}{\mathbb{Z}}
\newcommand{\rz}{\mathbb{R}}

\newcommand{\im}{\mathtt{i}}

\DeclareMathOperator{\Log}{Log}
\DeclareMathOperator{\Arg}{Arg}
\DeclareMathOperator{\lcm}{lcm}

\newtheorem{proclaim}{Lemma}

\def\udots{\mathinner{\mskip1mu\raise1pt\vbox{\kern7pt\hbox{.}}\mskip2mu
       \raise4pt\hbox{.}\mskip2mu\raise 7pt\hbox{.}\mskip1mu}}   

\begin{document}


\title{How to get high resolution results from sparse and coarsely sampled data}

\author{Annie Cuyt and Wen-shin Lee\\
Department of Mathematics and Computer Science\\ University of Antwerp (CMI)\\ 
Middelheimlaan 1, B-2020 Antwerpen, Belgium
\\ \texttt{\{annie.cuyt,wen-shin.lee\}@uantwerpen.be}}






\maketitle

\begin{abstract}
Sampling a signal below the Shannon-Nyquist rate causes aliasing, meaning
different frequencies to become indistinguishable. It is also well-known
that recovering spectral information from a signal using a parametric 
method can be ill-posed or ill-conditioned and therefore should be done
with caution. 

We present an exponential analysis method to retrieve high-resolution 
information from coarse-scale measurements, using uniform downsampling.
We exploit rather than avoid aliasing. While we loose the unicity of the
solution by the downsampling, it allows to recondition the problem
statement and increase the resolution. 

Our technique 
can be combined with different existing implementations of
multi-exponential analysis (matrix pencil, MUSIC, ESPRIT, APM, 
generalized overdetermined eigenvalue solver, simultaneous QR
factorization, $\ldots$) and so is very versatile.
It seems to be especially useful in the presence of clusters of
frequencies that are difficult to distinguish from one another.

\end{abstract}


\noindent\textbf{Mathematics Subject Classication (2010):} 42A15, 65Z05, 65T40.



\section{Introduction}\label{sect:one}
Estimating the fine scale spectral information of an exponential sum 
plays an important role in many signal processing applications. 
The problem of superresolution \cite{CandesFernandez:2014,Moitra:2015} 
has therefore recently received considerable attention.

Despite its computational
efficiency and wide applicability, the often used Fourier transform
(FT) has some well-known limitations,
such as its limited resolution and the leakage in the
frequency domain.
These restrictions complicate the analysis of 
signals falling exponentially with time. Fourier
analysis, which represents a signal as a sum of periodic functions, is
not very well suited for the decomposition of aperiodic signals, such as
exponentially decaying ones. The damping causes a broadening of the
spectral peaks, which in its turn leads to the peaks overlapping and
masking the smaller amplitude peaks. The latter are important for the
fine level signal classification. 

Signals that fall exponentially with time appear, for instance, in
transient detection,
motor fault diagnosis,
electrophysiology,
magnetic resonance and infrared spectroscopy,
vibration and
seismic data analysis,
music signal processing,
corrosion rate and crack initiation modelling,
electronic odour recognition, 
typed key\-stro\-ke recognition, 
nuclear science, liquid explosives identification, direction of arrival
estimation, and so on. 

A different approach to spectral analysis is offered, among others, 
by parametric methods.
However, parametric methods often suffer from ill-posedness 
and ill-conditioning 
particularly in the case of
clustered frequencies \cite{Kammler,Va:fit:85,Ha.Ka:eff:97}. 
In general, parametric methods
also require prior knowledge of the model order. Widely used parametric
methods assuming a multi-expo\-nential model include 
MUSIC \cite{Sc:mul:86}, ESPRIT \cite{Ro.Ka:esp:89}, the matrix pencil
algorithm \cite{Hu.Sa:mat:90}, simultaneous QR factorization
\cite{Go.Mi.ea:sta:99} or a generalized overdetermined eigenvalue 
solver \cite{Neumaier} and the approximate Prony method APM
\cite{Be.Mo:app:05,Po.Ta:par:10,Po.Ta:par:13}.

In general, parametric methods as well as the FT,
sample at a rate dictated by
the Shannon-Nyquist theorem \cite{Ny:cer:28,Sh:com:49}.
It states that the sampling rate needs to be
at least twice the maximum frequency present in the signal. A coarser time
grid causes aliasing,
identifying higher frequencies with lower frequencies without being able
to distinguish between them. Conventional
measurement systems, as used in modern consumer electronics, biomedical
monitoring and medical imaging devices, are all based on this theorem.

In the past decade, alternative approaches have proved that
signal reconstruction is also possible from sub-Nyquist measurements,
if additional
information on the structure of the signal is known, such as its
sparsity. 
Many signals are indeed spar\-se in some domain such as time, frequency or
space, meaning that most of the samples of either the signal or its
transform in another domain can be regarded as zero. Among others, we refer
to compressed sensing \cite{Ca.Ro.ea:rob:06,Donoho:2006}, finite rate of
innovation \cite{VMB02}, the use of coprime arrays in DOA
\cite{VaPa2011-2,TEN14}. 

The ultimate goal is to retrieve fine-scale information directly
from coarse-scale measurements acquired at a slower
information rate, in function of the sparsity and not the bandwidth
of the signal.
We offer a technique that allows to overcome the Shannon-Nyquist sampling
rate limitation and at the same time may improve the
conditioning of the numerical linear algebra problems involved. The
technique is exploiting aliasing rather than avoiding it and maintains a
regular sampling scheme \cite{Cu.Le:sma:12*b,Cu.Le:sma:12}. It relies on
the concept of what we call {\it identification shift}
\cite{Cu.Le:sma:12,Cu.Le:sma:12*b}, which is the additional sampling at
locations shifted with respect to the original locations, in order to
overcome any ambiguity in the analysis arising from periodicity issues
and in order to solve other identification
problems occurring in coprime array approaches. 

The paper is organized as follows. Exponential analysis following
Shannon-Nyquist sampling is repeated in Section~\ref{sect:two} and generalized to
sub-Nyquist sampling in Section~\ref{sect:three}. Since a sub-Nyquist rate can cause
terms to collide at the time of the sampling, we explain how to unravel
collisions in Section~\ref{sect:four}. Such collisions are very unlikely to happen
in practice of course.
In Section~\ref{sect:five} numerical examples illustrate both
the collision-free situation and the case in which the collision of terms
happens. The numerical examples at the same time illustrate:
\begin{itemize}
\item how the method
reconditions a problem statement that is numerically ill-conditioned
because of the presence of frequency clusters,
\item that it can be 
combined with an existing implementation of a multi-exponential spectral
analysis (we used ESPRIT \cite{Ro.Ka:esp:89} and {\tt oeig} \cite{Neumaier}).
\end{itemize}



\section{The multi-exponential model}\label{sect:two}
In order to proceed we introduce some notations. Let the real parameters
$\psi_i, \omega_i, \beta_i$ and $\gamma_i$ respectively denote
the damping, frequency, amplitude and phase in each component of the
signal
\begin{equation}
\phi(t) = \sum_{i=1}^n \alpha_i \exp(\phi_i t), \qquad \alpha_i =
\beta_i \exp(\im\gamma_i), \qquad \phi_i = \psi_i + \im 2\pi\omega_i.
\label{exp_mod}
\end{equation}
For the moment, we assume that the frequency content is limited by
\cite{Ny:cer:28,Sh:com:49}
$$|\Im(\phi_i)|/(2\pi)= |\omega_i| < \Omega/2, \qquad i=1, \ldots,n,$$
and we sample $\phi(t)$ at the equidistant points $t_j=j\Delta$ for
$j=0, 1, \ldots, 2n-1, \ldots$ with $\Delta \le 1/\Omega$. In the sequel
we
denote
$$f_j:= \phi(t_j), \qquad j=0, 1, \ldots, 2n-1, \ldots$$
The aim is to find the model order $n$, and the parameters
$\phi_1, \ldots, \phi_n$ and $\alpha_1, \ldots,
\alpha_n$ from the measurements $f_0, \ldots, f_{2n}, \ldots$ 
We further denote
$$\lambda_i:= \exp(\phi_i\Delta), \qquad i=1, \ldots, n.$$
With
$$H_n^{(k)} := \begin{pmatrix} f_k & \ldots & f_{k+n-1} \\ \vdots & \udots
& \vdots \\ f_{k+n-1} & \ldots & f_{k+2n-2} \end{pmatrix}, \qquad k \ge
0, \qquad n \ge 1,$$
the $\lambda_i$ are retrieved
\cite{Hu.Sa:mat:90,Go.Mi.ea:sta:99,Pl.Ta:pro:14} as the
generalized eigenvalues of the problem
\begin{equation}
H_n^{(1)} v_i = \lambda_i H_n^{(0)} v_i, \qquad i=1, \ldots,n,
\label{gep}
\end{equation}
where $v_i$ are the generalized right eigenvectors. From the values
$\lambda_i$, the complex numbers $\phi_i$ can be retrieved uniquely
because of the restriction $|\Im(\phi_i\Delta)|<\pi$.

In the absence of noise, the exact value for $n$ can be deduced from
\cite[p.~603]{He:app:74}, because we have for any single value of $k$ that
(for a detailed discussion also see \cite{Ka.Le:ear:03})
\begin{equation}
\begin{aligned}
\det H_\nu^{(k)} &=0 \text{ only accidentally}, \qquad k \ge 0, \\
\det H_n^{(k)} &\neq 0, \qquad k \ge 0 \\
\det H_\nu^{(k)} &= 0, \qquad \nu > n, \qquad k \ge 0.
\end{aligned} \label{det0}
\end{equation}
The way \eqref{det0} is checked is usually by computing the numerical rank
of a Hankel matrix $H_\nu^{(k)}$ or a rectangular $(\mu-\nu)\times \nu$ 
version of it with $\mu \ge 2\nu, \nu \ge n$, 
from its singular value decomposition \cite{Pl.Ta:pro:14}.
In the presence of noise and/or clusters of eigenvalues, this technique
may not be reliable though, 
but then some convergence property can be used instead \cite{Cu.Ts.ea:fai:18}. 
Note that hitting a zero value for $\det H_\nu^{(0)}$ accidentally, meaning
while $1\le\nu<n$, can only happen a finite number of times in a row, 
namely $n-1$ times (which is extremely unlikely), 
while the true value of $n$ is confirmed an infinite
number of times when overshooting it with any $\nu>n$. Therefore the
output of
\eqref{det0} is always probabilistic of nature. In Section 4.2 a similar
result is presented in the context of sub-Nyquist sampling where one may
loose the mutual distinctiveness of the generalized eigenvalues which is at the
basis of \eqref{det0}.

Finally, the $\alpha_i$ are computed from the interpolation conditions
\begin{equation}
\sum_{i=1}^n \alpha_i \exp(\phi_i t_j) = f_j, \qquad j=0, \ldots,
2n-1, \label{vdm}
\end{equation}
either by solving the system in the least squares sense, in the presence
of noise, or by solving a subset of $n$ (consecutive) interpolation
conditions in case of a noisefree $\phi(t)$. Also, $n$ can everywhere
be replaced by $N>n$, in order to model noise on the data by means of some
additional $N-n$ noise terms in \eqref{exp_mod}.
Note that
$$\exp(\phi_i t_j)=\lambda_i^j$$
and that the coefficient matrix of \eqref{vdm} is therefore a Vandermonde
matrix. It is well-known that the conditioning of structured matrices is
something that needs to be monitored \cite{Be.Go.ea:num:07,Ga:nor:75}.

Without loss of generality, we assume in the sequel that $0 \le \omega_i <
\Omega \in \nz, i=1, \ldots, n$ instead of $|\omega_i|<\Omega/2, i=1,
\ldots, n$. 
Also we assume in Section~\ref{sect:three} that $n$ is known or correctly detected as
indicated in \cite{Cu.Ts.ea:fai:18}. In Section~\ref{sect:four} we explain how to detect
$n$ concurrently with the computation of the $\phi_i$ and $\alpha_i$ from
sub-Nyquist data.

\section{Sub-Nyquist multi-exponential analysis}\label{sect:three}
Some basic result is first deduced for $n=1$. Afterwards this result is 
made use of for general $n$. The latter however, demands additional
developments.

\subsection{Dealing with a single frequency ($n=1$)}\label{sect:three.one}
At first we deal with some simple mathematical results, without caring
about computational issues.
When 
$$\phi(t) = \alpha \exp (\psi t + \im 2\pi \omega t ), \qquad 0 \le
\omega<\Omega,$$
and $\phi(t)$ is sampled at $t_j=0, \Delta, 2\Delta, \ldots$, with for
simplicity $\Delta=1/\Omega$,
then $\omega$ can uniquely be determined in $[0, \Omega)$ from the
samples. No periodicity problem occurs since $\omega\Delta < 1$ in the
generalized eigenvalue
$$\lambda=\exp(\psi\Delta)\exp(\im 2\pi\omega\Delta).$$
When $\phi(t)$ is sampled at multiples $t_{r_1j}=0, r_1\Delta, 2r_1\Delta,
\ldots$ with $1<r_1 \in \nz$, then there exist $r_1$ solutions for $\omega$ in
$[0, \Omega)$ since $0 \le 2\pi\omega r_1\Delta < 2r_1\pi$.
If $\phi(t)$ is also sampled at $t_{r_2j}=0, r_2\Delta, 2r_2\Delta, \ldots$
with $0<r_2\in\nz$, then one obtains another set containing $r_2$
solutions for $\omega$. Each solution set is extracted from the respective
generalized eigenvalues $\exp(\psi r_m\Delta)\exp(\im 2\pi\omega
r_m\Delta), m=1,2$ 
satisfying \eqref{gep} where
the first generalized eigenvalue problem is set up with the samples
$f_{r_1j} = \phi(0), \phi(r_1\Delta), \phi(2r_1\Delta), \ldots$ and the second
generalized eigenvalue problem with the samples 
$f_{r_2 j} = \phi(0), \phi(r_2\Delta), \phi(2r_2\Delta), \ldots$
In our write-up we have chosen not to add an index $r$ to the notation of
the Hankel matrices $H_n^{(0)}$ and $H_n^{(1)}$ when they are 
filled with samples taken at multiples
$t_{rj}=jr\Delta$ in order to not overload the
notation. From the context it is always clear which sequence of samples is
meant.

It is easy to show that, if in addition 
$\gcd(r_1, r _2)=1$, then $\omega$ is the unique
intersection of the two solution sets.
\begin{proclaim}\label{lem:one} 
Let $0 \le \omega < \Omega$ and $\Omega, r_1, r_2$ be nonzero positive
integers. If $\gcd(r_1, r_2)=1$ and if $\Delta = 1/\Omega$, then from the
values $\exp(\im 2\pi \omega r_1\Delta)$
and $\exp(\im 2 \pi\omega r_2\Delta)$, the frequency $\omega$ can uniquely
be recovered in $[0, \Omega)$.
\end{proclaim}
\begin{pf}
From the generalized eigenvalue $\lambda^{r_1}=\exp(\psi r_1\Delta)
\exp(\im 2 \pi r_1\Delta)$ we extract
$r_1$ solutions for $\omega$:
\begin{equation}
\omega = \omega^{(1)} + k {\Omega \over r_1}, \qquad 0 \le \omega^{(1)} <
{\Omega \over r_1}, \qquad k=0, \ldots, r_1-1. \label{r1}
\end{equation}
From the value $\lambda^{r_2}=\exp(\psi r_2\Delta)\exp(\im 2 \pi r_2\Delta)$ we extract
$r_2$ solutions for $\omega$:
\begin{equation}
\omega = \omega^{(2)} + \ell {\Omega \over r_2}, \qquad 0 \le \omega^{(2)} <
{\Omega \over r_2}, \qquad \ell=0, \ldots, r_2-1. \label{r2}
\end{equation}
Note that the frequency $\omega$ we are trying to identify satisfies both
\eqref{r1} and \eqref{r2}.
Remains to show that the common solution to \eqref{r1} and \eqref{r2} is
unique. Suppose we have two distinct values for $\omega$ that both 
satisfy \eqref{r1} and
\eqref{r2}. This implies that there exist two distinct $0\le k_1, k_2 <
r_1$ such that
\begin{equation}
\begin{aligned} \omega^{(1)}+k_1{\Omega \over r_1} &= \omega^{(2)} + \ell_1 {\Omega
\over r_2}, \\
\omega^{(1)}+k_2{\Omega \over r_1} &= \omega^{(2)} + \ell_2 {\Omega
\over r_2}, 
\end{aligned} \label{r12}
\end{equation}
with $0 \le \ell_1, \ell_2 < r_2$ and
$\ell_1 \not=\ell_2$ because $k_1\not=k_2$. From \eqref{r12} we deduce
$$k_1-k_2 = {(\ell_1-\ell_2) r_1 \over r_2} \not=0.$$
Hence $r_2$ divides $\ell_1-\ell_2$ because $\gcd(r_1,r_2)=1$. 
Since $\ell_1-\ell_2$ is bounded in absolute value by $r_2-1$, this is a
contradiction. \qed
\end{pf}

Furthermore, the element $\omega$ in the intersection can be obtained from 
the Euclidean algorithm.
\begin{proclaim}\label{lem:two}  
Let $0 \le \omega < \Omega$ and $\Omega, r_1, r_2$ be nonzero positive 
integers. If $\gcd(r_1, r_2)=1$ and if $\Delta = 1/\Omega$, then from the
values $\exp(\im 2\pi \omega r_1\Delta)$
and $\exp(\im 2 \pi\omega r_2\Delta)$, the frequency $\omega \in
[0, \Omega)$ is obtained as
\begin{equation}
\left( p_1 {\Log(\exp(\im 2\pi \omega r_1\Delta)) \over \im 2 \pi} +
p_2 {\Log(\exp(\im 2\pi \omega r_2\Delta)) \over \im 2 \pi} \right) \Omega
= \omega + h \Omega, \quad h \in \gz, \label{EA}
\end{equation}
where $p_1r_1 + p_2 r_2 = 1 \mod\Omega$ with $p_1, p_2 \in \gz$ and
$\Log(\cdot)$ denotes the principal branch of the complex logarithm.
\end{proclaim}
\begin{pf}
We use the same notation as in the proof of Lemma 1. So we have
\begin{align*}
\omega &= \omega^{(1)}+k{\Omega \over r_1}, \qquad 0 \le \omega^{(1)} <
{\Omega \over r_1}, \\
\omega &= \omega^{(2)}+\ell{\Omega \over r_2}, \qquad 0 \le \omega^{(2)} <
{\Omega \over r_2}.
\end{align*}
Then
$${\Omega \over r_1}
{\Log\left( \exp(\im 2\pi \omega r_1\Delta) \right) \over \im 2\pi} = 
\omega^{(1)}, \qquad
{\Omega \over r_2}
{\Log\left( \exp(\im 2\pi \omega r_2\Delta) \right) \over \im 2\pi} = 
\omega^{(2)},$$
and
\begin{multline*}
\left( p_1 {\Log\left( \exp(\im 2\pi \omega r_1\Delta) \right) \over \im
2\pi} +  p_2 {\Log\left( \exp(\im 2\pi \omega r_2\Delta) \right) \over \im
2\pi} \right) \Omega \\
\begin{aligned}
 &= (p_1r_1+p_2r_2)\omega - (p_1k+p_2\ell)\Omega \\
 &= \omega-(p_1k+p_2\ell)\Omega,
\end{aligned}
\end{multline*}
in which $p_1k+p_2\ell$ is an integer. \qed
\end{pf}

When the integers $p_1$ and $p_2$ are small this method is very useful.
Otherwise one has to be careful about the numerical stability of
\eqref{EA}. One can of
course experiment with different $r_1$ and $r_2$ values to ensure small
$p_1$ and $p_2$ values.

\subsection{Dealing with several terms ($n>1$)}\label{sect:three.two}
When $\phi(t)$ contains several terms, then we obtain $n$ solution sets
for the $\omega_i, i=1, \ldots, n$ from the first batch of evaluations at
multiples of $r_1\Delta$ and another $n$ solution sets for these
frequencies from the second batch of samples at multiples of $r_2\Delta$.
But now we are facing the problem of correctly matching the solution
set from the first batch to the solution set from the second batch 
that refer to the same $\omega_i$. Of course,
we want to avoid such combinatorial steps in our algorithm. 
To solve this problem we are going to choose the second batch of 
sampling points in a smarter way.

Before we proceed we assume that we don't have
$\exp(\phi_k r\Delta) = \exp(\phi_\ell r\Delta)$ for distinct $k$
and $\ell$ with $1 \le k, \ell \le n$. 
In Section~\ref{sect:four} we explain how to deal with the collision of terms, which we
exclude in the sequel of this section.

The sampling strategy that we propose is the following.
Sampling at $t_{rj}=jr\Delta$ with fixed $1<r\in\nz$,
gives us only aliased values for $\omega_i$, obtained from $\Log(\exp(\im
2\pi\omega r\Delta))$.
This aliasing can be fixed at the expense of the following 
additional samples.
In what follows $n$ can also everywhere be replaced by $N > n$ when using
$N-n$ additional terms in \eqref{exp_mod} to model the noise.

To fix the aliasing, we add $n$ samples to the already collected
$f_0, f_r, \ldots, f_{(2n-1)r}$,
namely at the shifted points
\begin{align*}
t_{rj+\rho} = j r\Delta + \rho & \Delta, \qquad r, \rho \text{ fixed,} \\
& j = h, \ldots, h+n-1, \qquad 0 \le h \le n.
\end{align*}
An easy choice for $\rho$ is a number mutually prime with $r$. 
For the most general choice allowed, we refer to \cite{Cu.Le:ana:17}. 
An easy practical generalization is when $r$ and $\rho$ are rational
numbers $r/s$ and $\rho/\sigma$ respectively with $r, s, \sigma \in \nz$
and $\rho \in \gz$. In that case the condition $\gcd(r, \rho)=1$ is
replaced by $\gcd(\overline r, \overline \rho)=1$ where $r/s=\overline
r/\tau, \rho/\sigma = \overline \rho / \tau$ with $\tau=\lcm(s, \sigma)$.
Also, the indices of the shifted points need
not be consecutive, but for ease of notation we assume this for now.

From the samples $f_0, f_r, \ldots, f_{(2n-1)r}$ we first compute the
generalized eigenvalues $\lambda_i^r = \exp{(\phi_i r \Delta)}$
and the coefficients $\alpha_i$ going with $\lambda_i^r$ in the model
\begin{equation}
\begin{aligned}
\phi(j r \Delta) &= \sum_{i=1}^n \alpha_i \exp(\phi_i jr\Delta) \\
&= \sum_{i=1}^n \alpha_i \lambda_i^{jr}, \quad j=0, \ldots, 2n-1.
\end{aligned} \label{model}
\end{equation}
So we know which coefficient $\alpha_i$ goes with which generalized
eigenvalue $\lambda_i^r$, but
we just cannot identify the correct $\Im(\phi_i)$ from $\lambda_i^r$. The
samples $f_{jr+\rho}$ at the additional points
$t_{rj+\rho}$ satisfy
\begin{equation}
\begin{aligned}
\phi(jr\Delta + \rho \Delta) &= \sum_{i=1}^n \alpha_i \exp \left( \phi_i
(jr+\rho) \Delta \right) \\
&= \sum_{i=1}^n (\alpha_i \lambda_i^{\rho}) \lambda_i^{jr}, \\
&\qquad \qquad j=h, \ldots, h+n-1, \qquad 0 \le h \le n, 
\end{aligned} \label{shift}
\end{equation}
which can be interpreted as a linear system with the same coefficient
matrix entries as (\ref{model}), but now with a new left hand
side and unknowns $\alpha_1 \lambda_1^{\rho}, \ldots, \alpha_n
\lambda_n^{\rho}$ instead of $\alpha_1, \ldots, \alpha_n$.
And again we can associate each computed
coefficient $\alpha_i \lambda_i^\rho$ with
the proper generalized eigenvalue $\lambda_i^r$.
Then by dividing the $\alpha_i \lambda_i^{\rho}$ computed from
(\ref{shift}) by the $\alpha_i$ computed from (\ref{model}), for $i=1,
\ldots, n$, we
obtain from $\lambda_i^\rho$ a second set of $\rho$ plausible values
for $\omega_i$. Because of the fact that we choose $\rho$ and $r$
relatively prime, the two sets of plausible values for $\omega_i$ have
only one value in their intersection, as explicited in Lemma 1 and 2. 
Thus the aliasing problem is solved. 

\section{When aliasing causes terms to collide}\label{sect:four}
When $\exp(\phi_k r \Delta) = \exp(\phi_\ell r \Delta)$ with $k \not=
\ell$, then different exponential terms in \eqref{model} collide into one term 
as a consequence of the undersampling and the aliasing effect. 
Note that then for the moduli of the exponential terms holds that
$\exp(\psi_k r\Delta)=\exp(\psi_\ell r\Delta)$
and consequently $\psi_k=\psi_\ell$. As long as $\psi_k
\not=\psi_\ell$, exponential terms can be distinguished on the basis of
their modulus. So our focus is on the situation where 
$$\phi_k = \psi_k + \im 2 \pi \omega_k \not= \phi_\ell=\psi_\ell + \im 2 \pi
\omega_\ell, \qquad \psi_k = \psi_\ell, \quad r\omega_k = r\omega_\ell + h
\Omega, \quad h \in \gz.$$
Since terms can collide when subsampling, their correct number $n$ may not
be revealed when sampling at multiples of $r\Delta$, in other words, when
sampling at the rate $\Omega/r$ instead of $\Omega$. Let us assume that
\eqref{det0}, or its practical implementation in \cite{Cu.Ts.ea:fai:18} on $N
\times N$ Hankel matrices with $N > n$, 
reveals a total of $n_0$ terms after the first batch of
evaluations at $t_{rj}=jr\Delta$ with fixed $r$. 
We call $\lambda_i^{(0)}$ the $n_0$ generalized eigenvalues of \eqref{gep}
computed from the $f_{jr}$ as in Section~\ref{sect:three}. Since some of the terms in
\eqref{model} may have collided, we have
\begin{equation}\phi(t_{rj}) = \sum_{i=1}^{n_0} \alpha_i^{(0)} \exp(\phi_i^{(0)}t_{rj})
\label{four.one}\end{equation}
with 
$$\lambda_i^{(0)}= \exp(\phi_i^{(0)} r \Delta), \qquad i=1, \ldots, n_0,$$
and some of the
$\alpha_i^{(0)}$ being sums of the $\alpha_i$ from \eqref{model}.
In Section~\ref{sect:four.one} we assume that all $\alpha_i^{(0)}$ are nonzero.
The case where some of the collisions have disappeared because of
cancellations in the coefficients, meaning that some of the
$\alpha_i^{(0)}, i=1, \ldots, n_0$ 
are zero, is dealt with in Section~\ref{sect:four.two}.

It should be clear that the value of $n_0$ depends on $r$, as can be seen
in the following simple example (nevertheless we do not want to burden the
notation $n_0$ with this evidence). Consider the function $\phi(t)$ given
by
$$\phi(t) = e^{\im 2\pi t} - e^{\im 2\pi 21t} + e^{\im 2\pi 41t} - e^{\im
2\pi 61t} + e^{\im 2\pi 11t} - e^{\im 2\pi 31t} + e^{\im 2\pi 51t}.$$
With $\Delta=1/100$ and $r=5$ we find that in the evaluations
$\phi(jr\Delta)$ the first 4 terms cancel each
other and the last 3 terms collide into
\begin{equation}
f_{5j} = e^{\im 2\pi 55j/100}, \qquad n_0=2. \label{nzero}
\end{equation}
With $\Delta=1/100$ and $r=12$ the fourth and the fifth term cancel each
other and the first and the last term collide, giving
$$f_{12j} = 2e^{\im 2\pi 12j/100} - e^{\im 2\pi 52j/100} + e^{\im 2\pi
92j/100} - e^{\im 2\pi 72j/100}, \qquad n_0=4.$$

\subsection{Collision without cancellation}\label{sect:four.one}
We remark that $n_0 \le n$ and that
the $\phi_i^{(0)}$ are definitely among the $n$ parameters $\phi_i$ 
in \eqref{model}. 
Without loss of generality we assume that the colliding terms are
successive,
\begin{multline*}
\begin{pmatrix} \alpha_1^{(0)} \\ 
\vdots \\
\alpha_{n_0}^{(0)} \end{pmatrix} = 
\begin{pmatrix} \alpha_{h_1} + \cdots + \alpha_{h_2-1} \\
\vdots \\
\alpha_{h_{n_0}} + \cdots + \alpha_{h_{n_0+1}-1} 
\end{pmatrix}, \\ h_1 = 1, \qquad h_i \le h_{i+1}, 
\qquad 1 \le i \le n_0, \qquad h_{n_0+1} = n+1. 
\end{multline*}
In brief, when collisions occur, the computations return the results 
\begin{equation}
\begin{aligned}
\alpha_i^{(0)} &= \sum_{\ell=h_i}^{h_{i+1}-1} \alpha_\ell, \qquad i=1,
\ldots, n_0 \\
\lambda_i^{(0)} &= \lambda_{h_i}^r = \ldots = \lambda_{h_{i+1}-1}^r,
\qquad i=1, \ldots, n_0.
\end{aligned} \label{four.two}
\end{equation} 
Note that only the nonzero $\alpha_i^{(0)}$ and the distinct
$\lambda_i^{(0)}$ are revealed, without any knowledge about the $h_i, 1
\le i\le n_0$. In Section~\ref{sect:four.two} 
we explain how to deal with the additional problem
where some of the $\alpha_i$ cancel each other and therefore some of the
$\lambda_i^{(0)}$ have gone missing in the samples $\phi(t_{rj})$. 

For the sake of completeness we explicit the linear system that delivered the
$\alpha_i^{(0)}$, namely
\begin{equation}
\begin{pmatrix} 1 & \ldots & 1 \\ \lambda_1^{(0)} & \ldots &
\lambda_{n_0}^{(0)} \\ \vdots & & \vdots \\
(\lambda_1^{(0)})^{n_0-1} & \ldots & 
(\lambda_{n_0}^{(0)})^{n_0-1}
\end{pmatrix}
\begin{pmatrix} \alpha_1^{(0)} \\ \vdots \\ \alpha_{n_0}^{(0)} \end{pmatrix}
= \begin{pmatrix} f_{0} \\ f_{r} \\ \vdots \\ f_{(n_0-1)r}
\end{pmatrix} \label{Vreuse}
\end{equation}
or, as is most often the case, an overdetermined version of it.
We now explain how to disentangle the collisions, again making use of some
additional samples at shifted locations. Let $r$ and $\rho$ be fixed as
before with $\gcd(r, \rho)=1$. If $\gcd(r, \rho)>1$ for some reason or
because of a practical constraint, then the procedure may be an iterative one, 
as we indicate further below. 

Let us sample $\phi(t)$
at the shifted locations $t_{rj+\rho k}=(jr+k\rho)\Delta,
j= 0, \ldots, n_0-1, k\ge 1$. These sample values equal
\begin{equation}
f_{jr+k\rho} := \sum_{i=1}^{n_0} \left(
\sum_{\ell=h_i}^{h_{i+1}-1} \alpha_\ell\exp(\phi_\ell k\rho\Delta) 
\right) \exp(\phi_i^{(0)}jr\Delta). \label{Sshift}
\end{equation}
In \eqref{Sshift} we abbreviate 
\begin{equation}
\alpha_{i}^{(1)}(k) := \sum_{\ell=h_i}^{h_{i+1}-1}
\alpha_\ell\exp(\phi_\ell k\rho\Delta), \qquad i=1, \ldots, n_0.
\label{alpha_i1}
\end{equation}
For $k=0$ we have $\alpha_i^{(1)}(0)=\alpha_i^{(0)}, i=1, \ldots, n_0$.
For fixed $k > 0$ the
values $\alpha_i^{(1)}(k), i=1, \ldots, n_0$ are obtained from
\eqref{Sshift} and
\begin{equation}
\begin{pmatrix} 1 & \ldots & 1 \\ \lambda_1^{(0)} & \ldots &
\lambda_{n_0}^{(0)} \\ \vdots & & \vdots \\
(\lambda_1^{(0)})^{n_0-1} & \ldots & 
(\lambda_{n_0}^{(0)})^{n_0-1}
\end{pmatrix} 
\begin{pmatrix} \alpha_1^{(1)}(k) \\ \vdots \\ \alpha_{n_0}^{(1)}(k) \end{pmatrix}
=\begin{pmatrix} f_{k\rho} \\ f_{r+k\rho} \\ \vdots \\ f_{(n_0-1)r+k\rho}
\end{pmatrix} \label{step1}
\end{equation}
or its least squares version.
The Vandermonde coefficient matrix of \eqref{step1} is the same as the one
used to compute $\alpha_i^{(0)}, i=1, \ldots, n_0$ in \eqref{Vreuse} 
from the samples $f_{jr}$, 
which is the case $k=0$. So the Vandermonde matrix is reused as
it is independent of the index $k$ appearing 
in the right hand side and in the vector of unknowns.

When collecting in this way, for each $1 \le i \le n_0$, the values
$\alpha_i^{(1)}(0)$, $\alpha_i^{(1)}(1)$, $\alpha_i^{(1)}(2)$, $\ldots$ we have
a separate exponential analysis problem per $i$, namely to identify the
number of terms in $\alpha_i^{(1)}(k)$ in \eqref{alpha_i1}. 
Note that the sampling rate used 
to collect the $\alpha_i^{(1)}(k)$ is $\Omega/\rho$. 

Now we fix $1 \le i \le n_0$ and proceed.
When the samples $\alpha_i^{(1)}(k)$ take the place of the values $f_k$ in
\eqref{gep} and \eqref{vdm} and $h_{i+1}-h_i$ that of $n$, 
then: 
\begin{itemize}
\item the generalized eigenvalue problem \eqref{gep} 
delivers the components $\lambda_\ell^{(1)} =
\exp(\phi_\ell \rho \Delta)$ in \eqref{alpha_i1}, 
\item and the respective Vandermonde
system \eqref{vdm} delivers the $\alpha_\ell$ for $\ell=h_i, \ldots,$
$h_{i+1}-1$. 
\end{itemize}
Both can again be set up in a least squares sense, in a similar way as
for the determination of the $\lambda_i^{(0)}$ and $\alpha_i^{(0)}$.
As shown in
Lemma 1, the exponential sums
$\alpha_i^{(1)}(k)$ are fully disentangled and all $n$ terms
in \eqref{exp_mod} are identified when $\gcd(r, \rho)=1$, which is what 
we try to achieve in practice. 

With 
\begin{equation}
\label{four.three}
\lambda_\ell^{(1)} = \exp(\phi_\ell \rho \Delta)=\lambda_\ell^\rho, 
\qquad \ell=h_i, \ldots, h_{i+1}-1, \quad i=1, \ldots, n_0,
\end{equation}
and
$$\lambda_i^{(0)} = \exp(\phi_{h_i} r \Delta) = \lambda_\ell^r, \qquad
\ell=h_i, \ldots, h_{i+1}-1,\quad i=1, \ldots, n_0,$$
we have what we need in order 
to identify the $\phi_i, i=1, \ldots, n$ using Lemma
2, since
$$n=\sum_{i=1}^{n_0} (h_{i+1}-h_i).$$ 
An illustration of the procedure
above is presented in Section~\ref{sect:five.two}.

When for one or other reason
$\gcd(r, \rho)=s \not=1$ then the above procedure needs to be repeated
with $r$ replaced by $s$ and $\rho$ replaced by a suitable $\sigma$.
Then again additional samples are collected at shifted locations $t_{sj+\sigma
k}=(js+k\sigma)\Delta$, namely 
$$f_{js+k\sigma} := \phi\left( (js+k\sigma)\Delta \right),$$
and the procedure is repeated from \eqref{Sshift} on.
When $\gcd(r, \rho, \sigma)=1$ the procedure ends, otherwise it continues
as described.

\subsection{Collision with cancellation}\label{sect:four.two}
To complete the method, we discuss the special situation where some of
the terms $\alpha_i\exp(\phi_i jr\Delta)$ 
cancel each other when evaluating \eqref{exp_mod} at the $t_{rj}$, a
situation which is illustrated in Section~\ref{sect:five.three}. 

So at the first batch of evaluations $f_{jr}$,
in addition to collision, one encounters cancellation for one or more
indices $i, 1 \le i \le n_0$, meaning that one or more $\alpha_i^{(0)} =
\alpha_i^{(1)}(0)=0$. The fundamental
question is whether the $\alpha_i^{(1)}(k)$ can
continue to evaluate to zero for all $k$ in the second shifted batch of
evaluations $f_{jr+k\rho}$ when $\gcd(r, \rho)=1$? The answer is no, 
not even when the $\phi_\ell$ in \eqref{alpha_i1} have the same decay rate, 
as becomes clear from the Lemmas 3 and 4 below. 

\begin{proclaim}\label{lem:three} 
Let for $\phi_k\not=\phi_\ell$ and $r\not=0$ hold that $\exp(\phi_k
r\Delta) = \exp(\phi_\ell r\Delta)$. If $\gcd(r, \rho)=1$ then 
$$\exp(\phi_k \rho\Delta) \not= \exp(\phi_\ell \rho\Delta).$$ 
\end{proclaim}
\begin{pf}
As pointed out it is sufficient to deal with the imaginary parts of
$\phi_k$ and $\phi_\ell$. We use a similar notation as in Lemma 1. The
proof is by contraposition.
From $\exp(\phi_k r\Delta) = \exp(\phi_\ell r\Delta)$ 
and $\exp(\phi_k \rho\Delta) = \exp(\phi_\ell \rho\Delta)$ 
we find that there
exist integers $p_k, p_\ell, q_k, q_\ell$ such that
\begin{align*}
\omega_k &= \omega^{(1)} + p_k {\Omega \over r}, \qquad 0 \le p_k \le r-1 \\
\omega_\ell &= \omega^{(1)} + p_\ell {\Omega \over r}, \qquad 0 \le p_\ell \le r-1 \\
\omega_k &= \omega^{(2)} + q_k {\Omega \over \rho}, \qquad 0 \le p_k \le \rho-1 \\
\omega_\ell &= \omega^{(2)} + q_\ell {\Omega \over \rho}, \qquad 0 \le p_k
\le \rho-1 .
\end{align*}
Then
$$\omega_k-\omega_\ell = (p_k-p_\ell) {\Omega \over r} = (q_k-q_\ell)
{\Omega \over \rho}$$
or
$$p_k-p_\ell = {q_k-q_\ell \over \rho} r,$$
which is a contradiction since the left hand side is an integer and
$q_k-q_\ell$ in the right hand side is  in absolute value bounded by
$\rho-1$. \qed
\end{pf}
\begin{proclaim}\label{lem:four} 
Let $\alpha_i^{(1)}(k)$ be
given by \eqref{alpha_i1}. Then $\forall\; 1 \le i \le n_0, \exists\; 
0 \le k \le h_{i+1}-h_i:  \alpha_i^{(1)}(k) \not=0$.
\end{proclaim}
\begin{pf}
We consider the following
square Vandermonde system which is obtained from \eqref{alpha_i1} 
for fixed $i$ and with $k$ increased from 0 to $h_{i+1}-h_i$,
\begin{multline}
\label{qedemo}
\begin{pmatrix} 1 & \ldots & 1 \\ \exp(\phi_{h_i} \rho \Delta) & \ldots &
\exp(\phi_{h_{i+1}-1} \rho \Delta) \\ \vdots & & \vdots \\
\exp(\phi_{h_i} (h_{i+1}-h_i)\rho \Delta) & \ldots &
\exp(\phi_{h_{i+1}-1} (h_{i+1}-h_i)\rho \Delta) 
\end{pmatrix} 
\begin{pmatrix} \alpha_{h_i} \\ \vdots \\ \alpha_{h_{i+1}-1} \end{pmatrix} = \\
\begin{pmatrix} \alpha_i^{(1)}(0) \\ \vdots \\ \alpha_i^{(1)}(h_{i+1}-h_i)
\end{pmatrix}.
\end{multline} 
From Lemma 3 we know that this 
$(h_{i+1}-h_i) \times (h_{i+1}-h_i)$ Vandermonde matrix is regular.
If the right hand side of this small linear system
consists of all zeroes, then we must therefore conclude incorrectly that
$$\alpha_{h_i} = \cdots = \alpha_{h_{i+1}-1} = 0.$$
So from this we know that the evaluation of $\alpha_i^{(1)}(k)$
to zero cannot persist up to and including $k=h_{i+1}-h_i$. \qed
\end{pf}

The important conclusion here is that in a finite number of steps the true value
of $n_0$, which represents the number of distinct generalized eigenvalues
existing at the sampling rate $r\Delta$, is always revealed. The 
evaluations at the shifted sample points $(jr+k\rho)\Delta$ with $k \not= 0$ 
serve the purpose to  
provide a different view on the coeficients, namely the values
$\alpha_i^{(1)}(k)$ for $k \not=0$. These additional evaluations do not
alter or touch the generalized eigenvalues $\lambda_i^{(0)}$. 
That is why a shift is so
helpful. And Lemma 4 confirms, that even if initially some
$\alpha_i^{(1)}(0)$ are zero, eventually all $\alpha_i^{(1)}(k)$
must become visible. This fact is entirely similar to the conclusion in
\eqref{det0}, but with the function $\phi(t)$ replaced by 
$\alpha_i^{(1)}(k)$ for some fixed $i$ and with the matrix $H^{(0)}_\nu$ 
replaced by the matrix
$$\begin{pmatrix} \alpha_i^{(1)}(0) & \ldots & \alpha_i^{(1)}(\kappa-1) \\
\vdots & & \vdots \\
\alpha_i^{(1)}(\kappa-1) & \ldots & \alpha_i^{(1)}(2\kappa-2) 
\end{pmatrix}$$
of increasing size $\kappa\times\kappa$. 

To illustrate this we return to \eqref{nzero}. While only one of the
$n_0=2$ terms is visible when evaluating at $jr\Delta$ when $r=5$, the
evaluations $f_{5j+12k}$ with $\rho=12$, give us
\begin{align*}
f_{5j+12k} &= \left( e^{\im 2\pi 12k/100} - e^{\im 2\pi 52k/100} + e^{\im
2\pi 92k/100} - e^{\im 2\pi 32k/100} \right) e^{\im 2\pi 5j/100} + \\
& \hskip 1truecm 
\left( e^{\im 2\pi 32k/100} - e^{\im 2\pi 72k/100} + e^{\im 2\pi
12k/100} \right) e^{\im 2\pi 55j/100} \\
&= \alpha_{h_1}^{(1)}(k) e^{\im 2\pi 5j/100} + \alpha_{h_2}^{(1)}(k)
e^{\im 2\pi 55j/100}. 
\end{align*}
For $k \ge 1$ and $\nu \ge 2$ we find that the rank of 
$$H_\nu^{(12k)} \begin{pmatrix} f_{12k} & f_{r+12k} & \ldots & 
f_{(\nu-1)r+12k} \\ f_{r+12k} & & & \\ \vdots & & & \vdots \\
f_{(\nu-1)r+12k} & \ldots & & \end{pmatrix}$$ 
equals $n_0=2$.

\section{Numerical illustration}\label{sect:five}
We illustrate the working of \eqref{model} and \eqref{shift} from Section~\ref{sect:three} and
that of \eqref{Vreuse} and \eqref{step1} from Section~\ref{sect:four} on
two examples in the respective Sections~\ref{sect:five.one} and 
\ref{sect:five.two}. 
In the former numerical example the undersampling will not cause collisions, 
while in the latter illustration it will. In addition, in Section 
\ref{sect:five.three}, we show the detection of
terms that have not only collided but entirely vanished in the first sampling 
at multiples of $r\Delta$. We conclude in Section~\ref{sect:five.four} 
with pseudocode for
the full-blown algorithm, which is most easy to understand after going through
the examples. The pseudocode deals with all possible combinations of
situations and is therefore even more general than the example in Section
\ref{sect:five.three}.

\subsection{Collision-free example}\label{sect:five.one}
\begin{figure}
\centering\includegraphics[width=10truecm]{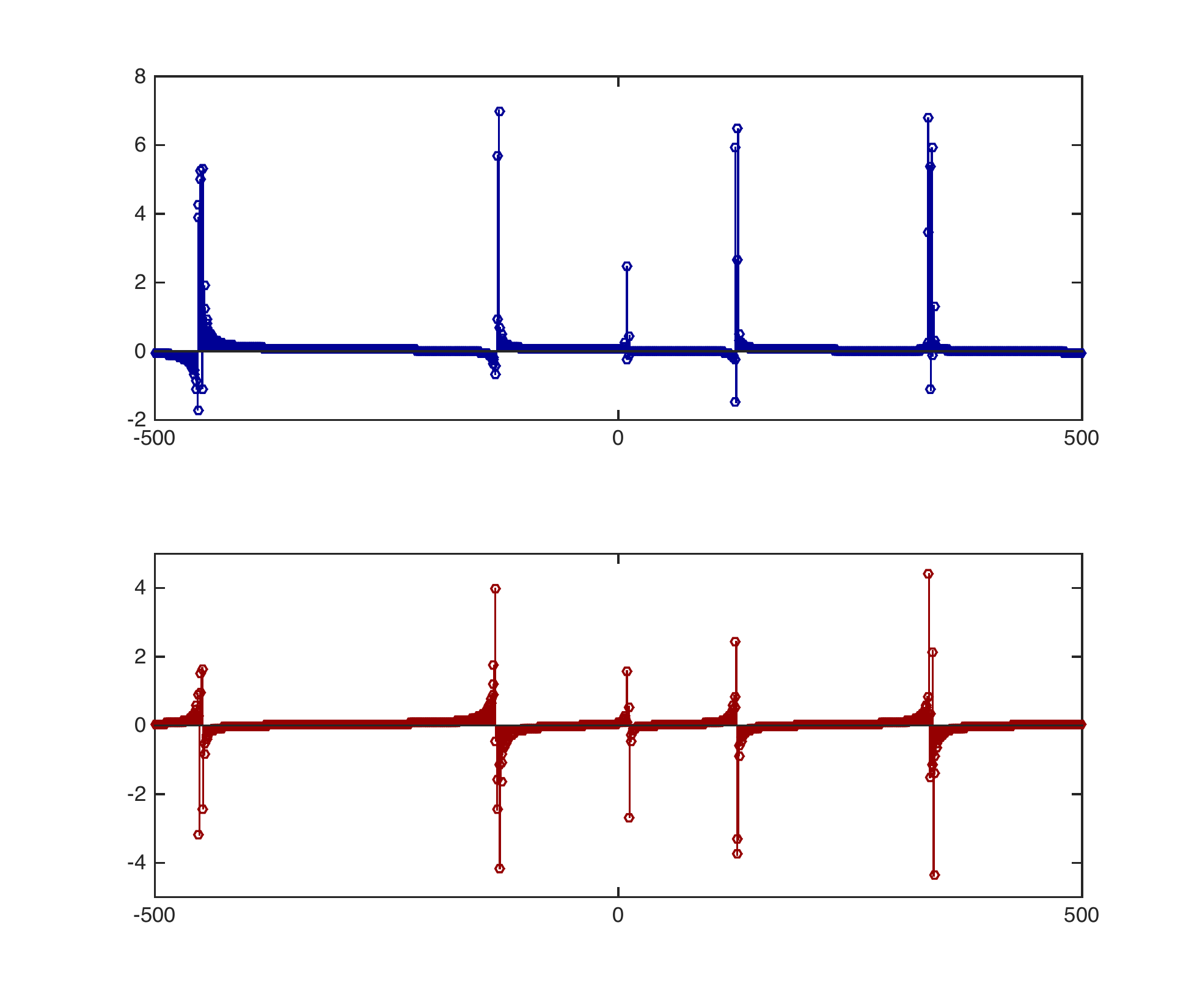}
\caption{Real (blue) and imaginary (red) 
part of the FT.}\label{fig:one}
\end{figure}

\begin{figure}
\centering\includegraphics[width=10truecm]{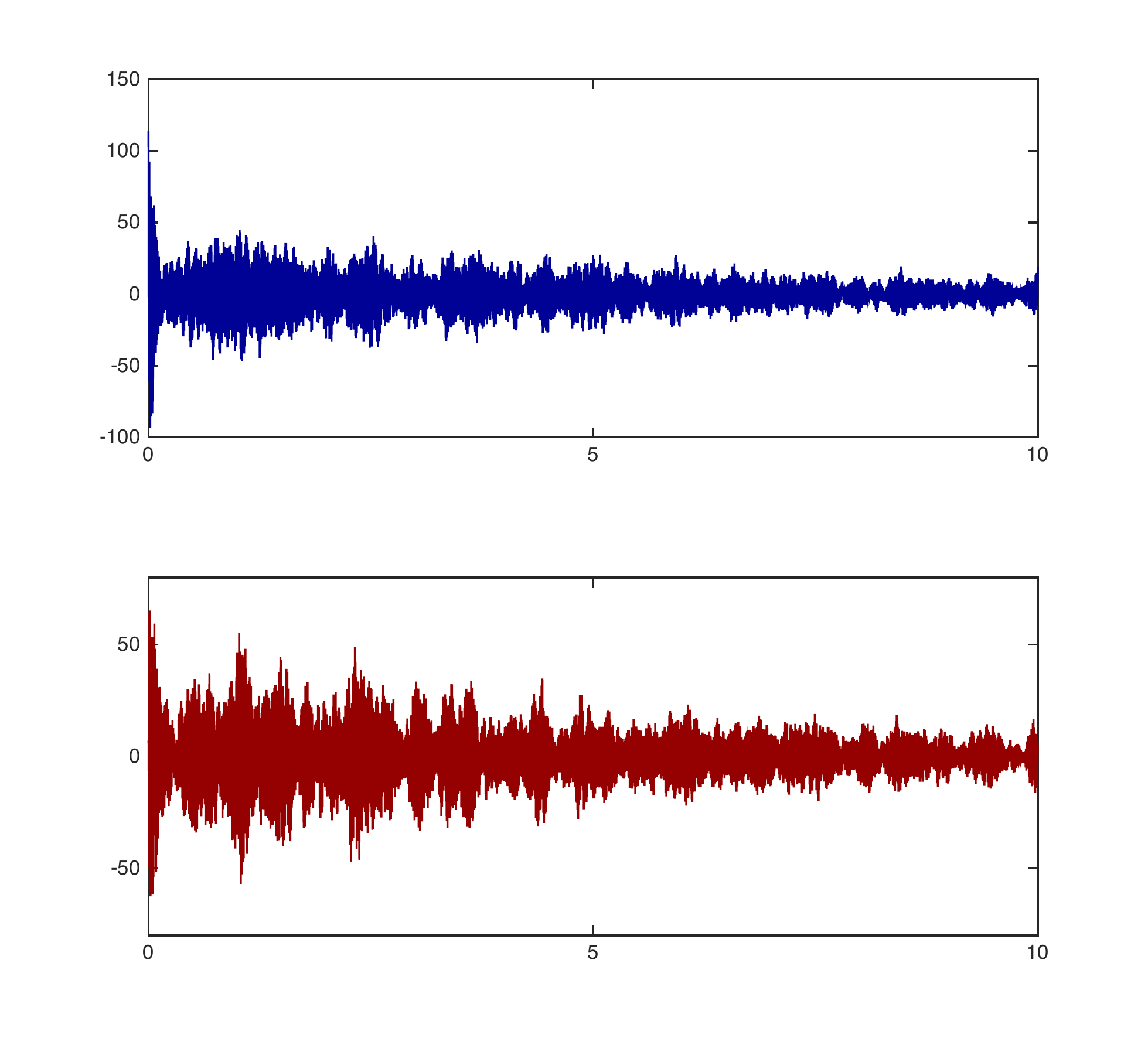}
\caption{Real (blue) and imaginary (red) part of the
signal.}\label{fig:two}
\end{figure}

\begin{figure}
\centering\includegraphics[width=9.5truecm]{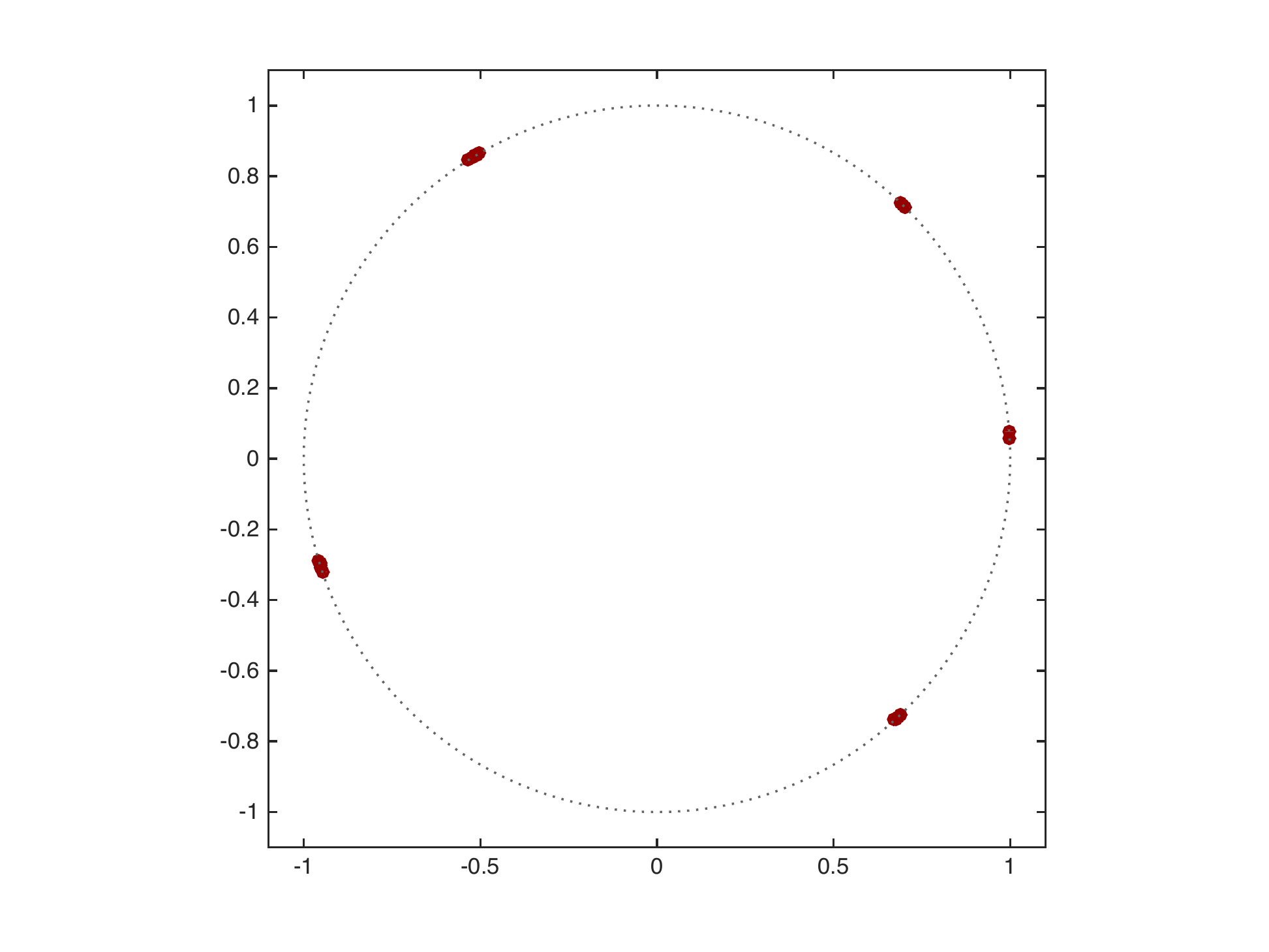}
\caption{Generalized eigenvalues $\lambda_i$ for
\eqref{exp_mod} with data from Table~\ref{table:one}.}\label{fig:three}
\end{figure}

For our first example the $\alpha_i$ and $\phi_i$ are given in Table~\ref{table:one}.
We take $\Omega=1000$ and $\Delta=1/\Omega$.
The $n=20$ frequencies $\omega_i$ form 5 clusters, as is apparent
from the FT, computed from 1000 samples and
shown in Figure~\ref{fig:one}. For completeness we graph the signal in
Figure~\ref{fig:two}. 
In Figure~\ref{fig:three} we show the generalized eigenvalues
$\lambda_i=\exp(\phi_i\Delta), i=1, \ldots,
20$ computed from the noisefree samples, to illustrate the
ill-conditioning of the problem as a result of the clustering of the
frequencies. 

\begin{figure}
\centering\includegraphics[width=10truecm]{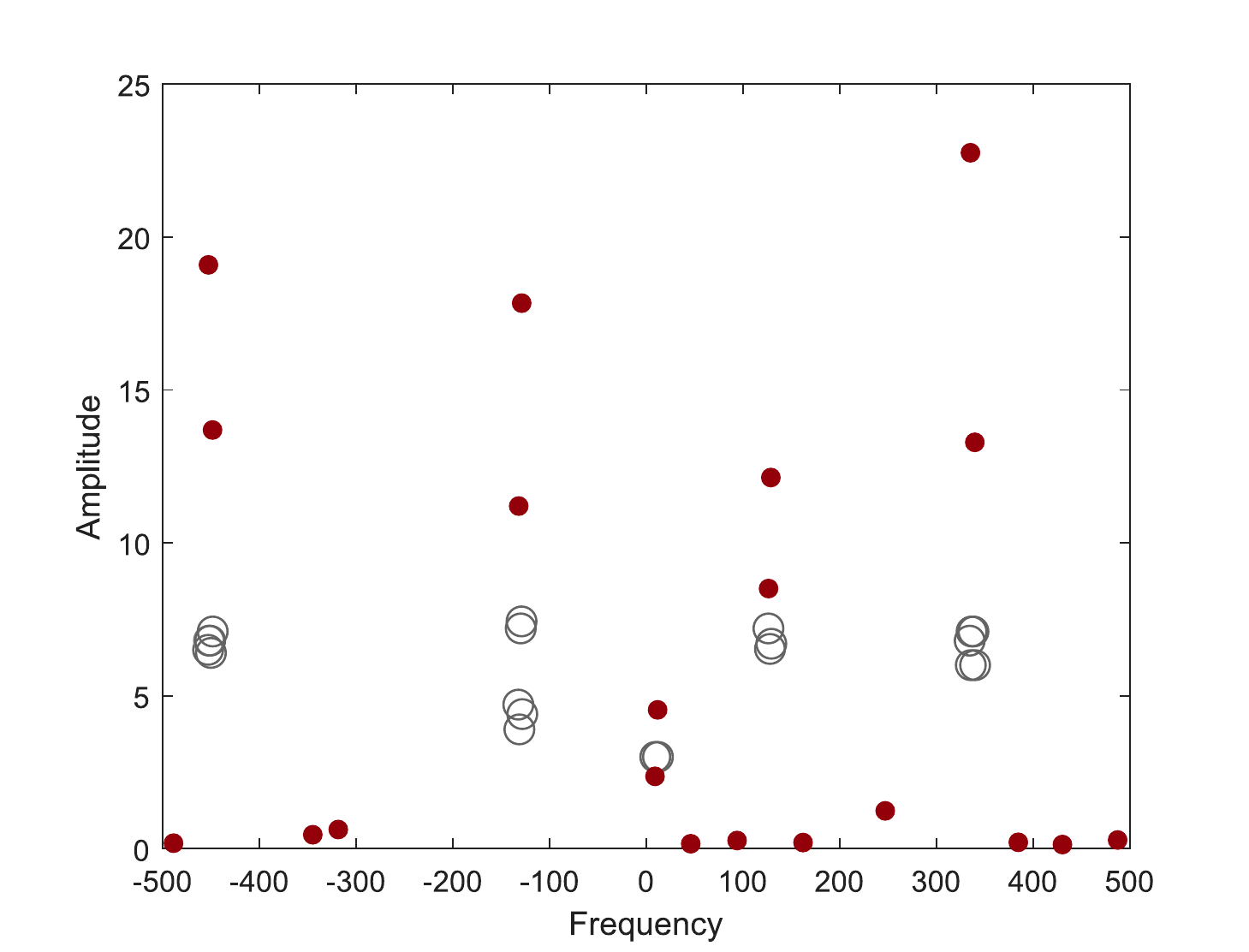}
\caption{ESPRIT output $(\omega_i, \beta_i), i=1, \ldots,
20$ computed from $f_0, \ldots, f_{239}$.}\label{fig:four}
\end{figure}

To \eqref{exp_mod} we add white Gaussian noise with SNR$=32$ dB. 
For comparison with our new method,
we show in Figure~\ref{fig:four} the $(\omega_i, \beta_i)$ results computed by means of 
ESPRIT using 240 samples, namely $f_0, \ldots, f_{239}$. A signal space
of dimension 20 and a noise space of dimension 40, so a total dimension
$N=60$, produced a typical ESPRIT result, from a 
size $180 \times 60$ problem. The true $(\omega_i,
\beta_i)$ couples from Table~\ref{table:one} are indicated using black circles. 
The ESPRIT output is indicated using red bullets. So ideally every black
circle should be hit by a red bullet. The ill-conditioning has
clearly created a serious problem in identifying the individual
input frequencies and amplitudes.

\begin{figure}
\centering\includegraphics[width=9.5truecm]{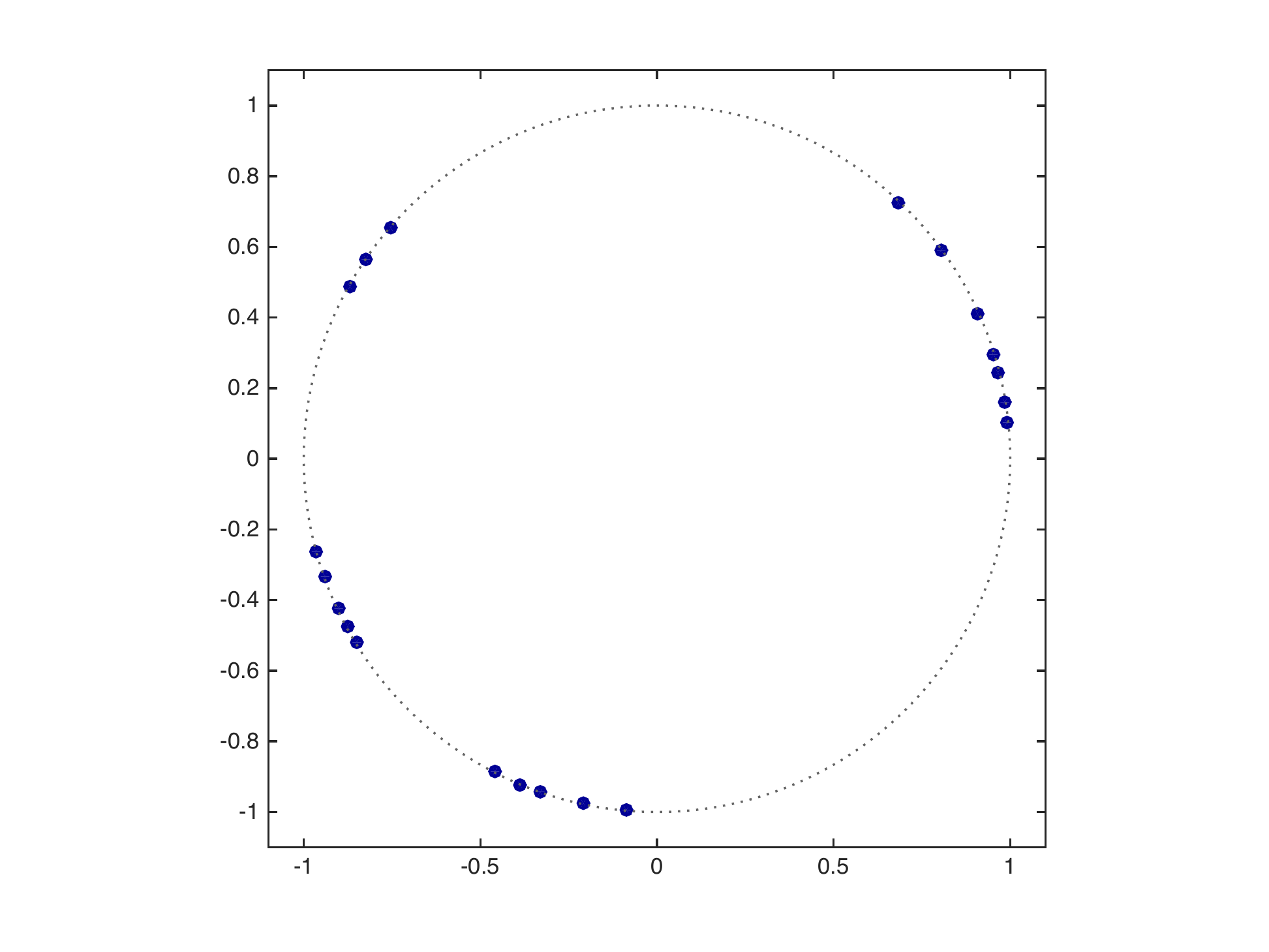}
\caption{Generalized eigenvalues $\lambda_i^r$ for
\eqref{exp_mod} with data from Table~\ref{table:one} and $r=11$.}\label{fig:five}
\end{figure}

Next we choose $r=11$ and $\rho=5$. The originally clustered eigenvalues
are now much better separated. To illustrate this we show 
in Figure~\ref{fig:five} the noisefree generalized
eigenvalues $\lambda_i^r, i=1, \ldots, 20$ of the $r$-fold undersampled 
exponential analysis problem. 

\begin{figure}
\centering\includegraphics[width=10truecm]{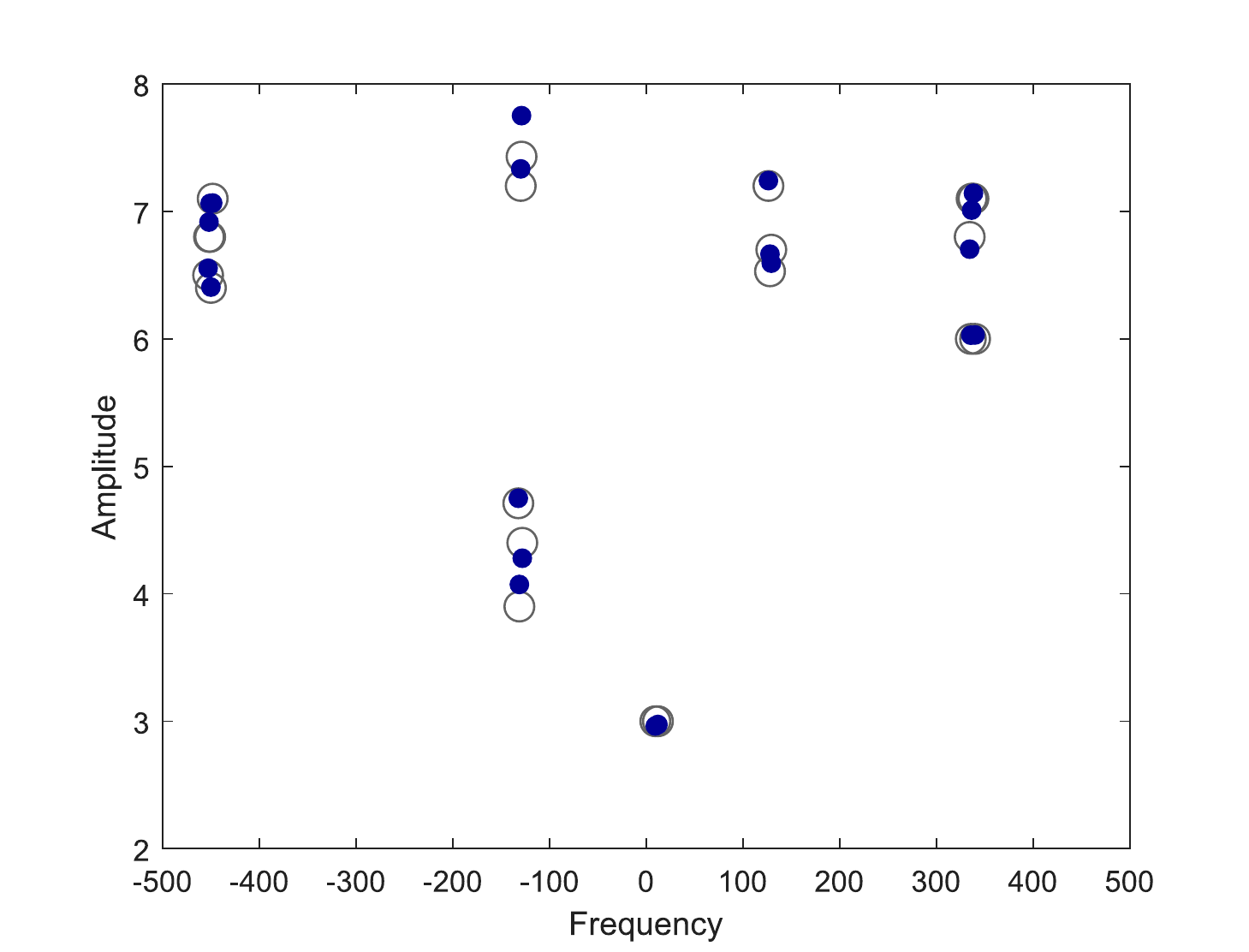}
\caption{Output $(\omega_i, \beta_i), i=1, \ldots,
20$ computed from 240 samples $f_{jr+\rho}, r=11, \rho=5$.}\label{fig:six}
\end{figure}

With the noisy samples, we again
take $N=60> n=20$ and set up a $120 \times 60$
generalized eigenvalue problem \eqref{gep} with the samples $f_{jr}, j=0,
\ldots, 179$ and the $120 \times 60$
Vandermonde system \eqref{model} that respectively deliver the
$\lambda_i^r$ and the $\alpha_i$ for $i=1, \ldots, N$. 
With the samples $f_{jr+\rho}, j=0,
\ldots, 59$ we set up the $60 \times 60$ linear system \eqref{shift} 
from which we compute the $\alpha_i \lambda_i^\rho, i=1, \ldots, 60$
and subsequently the $\lambda_i^\rho$. This brings our total number of
samples used also to 240, comparable to the ESPRIT procedure. An advantage for
ESPRIT is that the signal has less decayed in the first 240 samples,
compared to the 240 samples used here.
Using the Euclidean algorithm, as explicited in Lemma 2, we recover from
$\lambda_i^r$ and $\lambda_i^\rho$ the
true frequencies $\omega_i$ with $p_1=1, p_2=-2$ and 
$p_1 r + p_2\rho =1$. With the new method we find the $(\omega_i,
\beta_i)$ couples shown as blue dots in Figure~\ref{fig:six}. 
In Figure~\ref{fig:six} the reader can even
clearly count the number of frequencies retrieved in each cluster, which
is the correct number when comparing to the input values in Table~\ref{table:one}.
Clearly Figure~\ref{fig:six} is a tremendous improvement over Figure~\ref{fig:four}. 

\subsection{Example where collisions occur without cancellation}
\label{sect:five.two}
In Table~\ref{table:two} we list the $\alpha_i$ and $\phi_i$ of an exponential model,
chosen in such a way that the aliasing 
causes terms to collide. This enables us to illustrate the workings of the
technique explained in Section~\ref{sect:four}. 

\begin{figure}
\centering\includegraphics[width=7truecm]{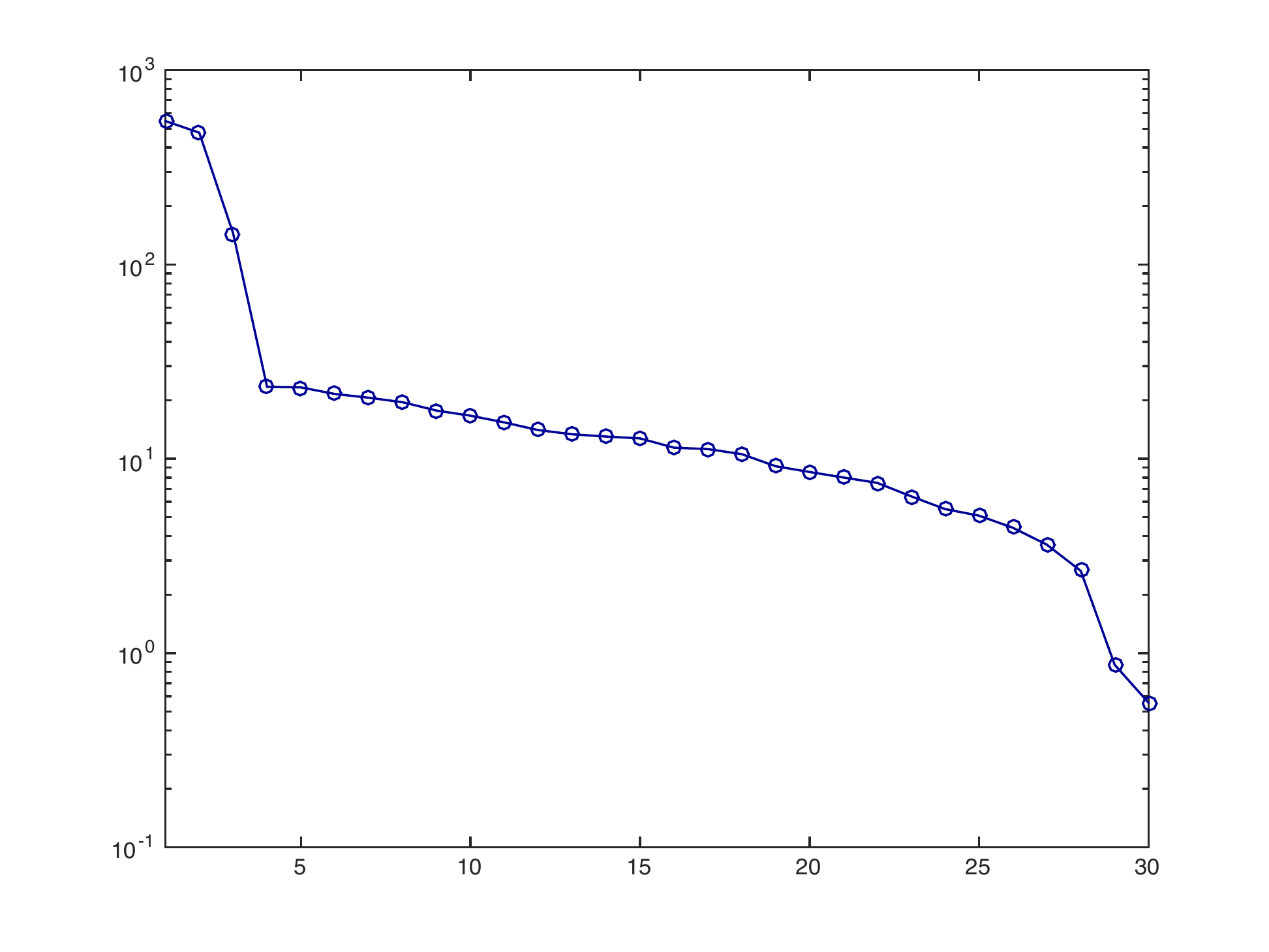}
\caption{SVD of $H_{30}^{(0)}$
for \eqref{exp_mod} with data from Table~\ref{table:two} and $r=100$.}\label{fig:seven}
\end{figure}

The bandwidth is again $\Omega=1000$ and we take $\Delta=1/\Omega$ and
$r=100$. 
We add white Gaussian noise to the samples with SNR$=20$ dB and start our
computations. When subsampling, the 6 terms collide into 3, as indicated 
in Figure~\ref{fig:seven} by the
singular value decomposition of $H^{(0)}_N$ with $N=30$,
which reveals its numerical rank. Actually
$$\phi(t_{rj}) = (\alpha_1+\alpha_2+\alpha_3) \exp(\phi_1 jr\Delta) + 
\alpha_4 \exp(\phi_4 jr\Delta) + (\alpha_5+\alpha_6)\exp(\phi_5
jr\Delta).$$
We recall that $H_N^{(0)}$ is filled
with the samples $f_{jr}, j=0, \ldots, 59$ and not with the samples $f_j,
j=0, \ldots, 59$.

We set up the $30 \times 30$ generalized
eigenvalue problem \eqref{gep} with the samples $f_{jr}, j=0, \ldots, 59$
which we solve using {\tt oeig},
and the $60 \times 30$ Vandermonde system \eqref{model} that respectively
deliver the $\lambda_i^r$ and the $\alpha_i$ for $i=1, \ldots, N$. When
retaining the components with largest $|\alpha_i|$, we find
\begin{align*}
\lambda_1^{(0)} &\approx 0.36845+0.93042 \im \\
\lambda_2^{(0)} &\approx 0.36745-0.92977 \im \\
\lambda_3^{(0)} &\approx -0.72761-0.68801 \im 
\end{align*}
and
\begin{equation}
\begin{aligned}
\alpha_1^{(0)} =\alpha_1+\alpha_2+\alpha_3 &\approx 17.718+0.25273 \im \\
\alpha_2^{(0)} =\alpha_4 &\approx 16.126+0.057118 \im \\
\alpha_3^{(0)} =\alpha_5 + \alpha_6 &\approx 4.5732-0.53331 \im
\end{aligned} \label{ssttaarr}
\end{equation}
At this point we have not yet been able to recover the correct $\lambda_i$ and
$\alpha_i$ for the signal defined by the parameters in Table~\ref{table:two} (we have
unearthed only 3 terms instead of 6) because of
two reasons. First, the subsampling creates an aliasing effect and second
the aliasing causes frequencies to collide. As explained in Section~\ref{sect:four}, we
can disentangle the information in the collisions from more values
$\alpha_i^{(1)}(k), k=1, 2, \ldots$, where 
$\alpha_i^{(1)}(0) = \alpha_i^{(0)}$, simply because the $\alpha_i^{(1)}(k)$
are themselves linear combinations of exponentials. 
To not complicate matters too much yet, the example is cancellation
free: so the correct value $n_0=3$ is immediately discovered from the
sampling at the multiples of $r\Delta$, as we see in \eqref{ssttaarr}.

For the disentanglement, we choose $\rho=133$ and
we set up the Vandermonde systems \eqref{step1},
$$\begin{pmatrix} 1 & \ldots & 1 \\ \lambda_1^{(0)} & \ldots &
\lambda_3^{(0)} \\ \vdots & & \vdots \\ (\lambda_1^{(0)})^9 & \ldots & 
(\lambda_3^{(0)})^9
\end{pmatrix} \begin{pmatrix} \alpha_1^{(1)}(k) \\ \alpha_2^{(1)}(k)
\\
\alpha_3^{(1)}(k) \end{pmatrix} = \begin{pmatrix} f_{k\rho} \\ f_{r + k\rho} \\ 
\vdots \\ f_{9r + k\rho} \end{pmatrix}, \qquad k=1, \ldots, 11.$$
In total so far 170 samples are used.
A singular value analysis of the Hankel matrices
$$\begin{pmatrix}
\alpha_i^{(1)}(0) & \alpha_i^{(1)}(1) & \ldots & \alpha_i^{(1)}(5) \\
\alpha_i^{(1)}(1) & \alpha_i^{(1)}(2) & \ldots & \alpha_i^{(1)}(6) \\
\vdots & \vdots & \ddots & \vdots \\
\alpha_i^{(1)}(5) & \alpha_i^{(1)}(6) & \ldots & \alpha_i^{(1)}(10) \\
\end{pmatrix}, \qquad i=1, 2, 3$$
reveals the number of components that one can distinguish in and
consequently extract from the $\alpha_i^{(1)}(k)$. 
The numbers are respectively 3, 2, 1 for $i=1, 2, 3$ and so $h_1=1, h_2=4,
h_3=6, h_4=7$. The size of these Hankel matrices filled with values
$\alpha_i^{(1)}(k)$, is chosen somewhat larger than necessary so that the
correctness of their rank is confirmed a number of times. We can also
conclude that
$$n = \sum_{i=1}^{n_0} (h_{i+1}-h_i) = 6.$$
For $i=1, 2, 3$ the generalized eigenvalue problems
$$\begin{pmatrix}
\alpha_i^{(1)}(1) &  \ldots & \alpha_i^{(1)}(6) \\
\vdots & \ddots & \vdots \\
\alpha_i^{(1)}(6) &  \ldots & \alpha_i^{(1)}(11) \\
\end{pmatrix} v_\ell
= \lambda_\ell^{(1)} \begin{pmatrix}
\alpha_i^{(1)}(0) & \ldots & \alpha_i^{(1)}(5) \\
\vdots  & \ddots & \vdots \\
\alpha_i^{(1)}(5) & \ldots & \alpha_i^{(1)}(10) \\
\end{pmatrix} v_\ell$$
reveal the $\lambda_\ell^{(1)}= \exp(\phi_\ell \rho \Delta), \ell=h_i,
\ldots, h_{i+1}-1, i=1, \ldots, n_0$. Note that we chose a notation where
the $\lambda_\ell^{(1)}$ are
not indexed by a double index $(\ell,i), \ell=1, \ldots, h_{i+1}-h_i, i=1,
\ldots, n_0$ but are indexed consecutively from $\ell=h_1=1$ to
$\ell=h_{n_0+1}-1=n$. This matches the indexing of the $\lambda_\ell^{(0)}$ 
of which some are coalescent, namely $\lambda_{h_i}^{(0)} = \cdots =
\lambda_{h_{i+1}-1}^{(0)}, i=1, \ldots, n_0$. 
The respective Vandermonde systems with unknowns
$\alpha_{h_i}, \ldots, \alpha_{h_{i+1}-1}$ and right hand sides
$\alpha_i^{(1)}(0)$, $\ldots$, 
$\alpha_i^{(1)}(11)$ reveal the $\alpha_\ell, \ell=h_i, \ldots, h_{i+1}-1$ in \eqref{alpha_i1}. 
Again we retain only the $h_{i+1}-h_i$ components with largest
$|\alpha_\ell|$.
From the $\lambda_\ell^{(0)} = \exp(\phi_\ell r \Delta), \ell=h_i,
\ldots, h_{i+1}-1, i=1, \ldots, n_0$ 
and $\lambda_\ell^{(1)}
= \exp(\phi_\ell \rho \Delta), \ell=1, \ldots, n$ the imaginary part 
of $\phi_i$ can be
recovered as indicated in Lemma 2: with $p_1=4$ and $p_2=-3$ we have $p_1
r + p_2 \rho = 1$ and so
$$\Im(\phi_\ell) =  4 \Arg(\lambda_\ell^{(0)})\Omega 
- 3 \Arg(\lambda_\ell^{(1)})\Omega + 2\pi h \Omega, \qquad 1 \le \ell \le
6, \qquad h \in \gz,$$
where $h$ is taken such that $0 \le \Im(\phi_\ell) < 2 \pi\Omega$.
Eventually we unearth the following 6 $\phi_i$ and $\alpha_i$:
\begin{align*}
\phi_1 &\approx -0.021600+ \im 2\pi 192.29, \\
\phi_2 &\approx -0.0085122+ \im 2\pi 289.87, \\
\phi_3 &\approx -0.025728+ \im 2\pi 386.69, \\
\phi_4 &\approx -0.066292+ \im 2\pi 538.18, \\
\phi_5 &\approx -0.043745+ \im 2\pi 858.70, \\
\phi_6 &\approx 0.0026126+ \im 2\pi 956.23 
\end{align*}
and
\begin{align*}
\alpha_1 &\approx 19.011 + \im 0.53818, \\
\alpha_2 &\approx -20.481 + \im 0.89352, \\
\alpha_3 &\approx 21.445 - \im 1.5790, \\
\alpha_4 &\approx 5.8439 - \im 0.035907, \\
\alpha_5 &\approx 5.0770 + \im 0.037562, \\
\alpha_6 &\approx 10.758 - \im 0.40878. 
\end{align*}

\subsection{Example with cancellations in the collisions}
\label{sect:five.three}

The cancellation strategy is most clearly illustrated 
by means of a noisefree example,
where exact cancellations are observed. The actual occurrence of this situation
in case of real-life data is extremely small, but we primarily want to show
that the proposed sub-Nyquist method is capable of recovering from it.

Let
\begin{multline}
\phi(t) = \exp(2\pi \im t) - \exp(2\pi\im 21 t) + \exp(2\pi\im
41 t ) - \exp(2\pi\im 61 t) +  \\ e^{\im 2\pi 72/100}\exp(2\pi\im 11t) 
- e^{\im 2 \pi 32/100} \exp(2\pi\im 31 t)+ \exp(2\pi\im 9 t). 
\label{excancel}
\end{multline}
We take $\Omega=100, \Delta=0.01$ and sample $f_j=\phi(j\Delta)$ for
particular values of $j$. With
$r=5$ the first four terms cancel each other and the fifth and sixth term
collide:
\begin{multline}
f_{5j} = 0 \exp(2\pi\im 5j/100) + \left( e^{\im 2\pi 72/100} - e^{\im
2\pi 32/100} \right) \exp(2\pi\im 55j/100) + \\ \exp(2\pi\im 45j/100).
\end{multline} 
So from the samples 
$f_0, f_5, f_{10}, \ldots$ only two terms can be retrieved:
$$\text{rank } H^{(0)}_3 = \text{rank} \begin{pmatrix} f_0 & f_5 & f_{10} \\
f_5 & f_{10} & f_{15} \\ f_{10} & f_{15} & f_{20}
\end{pmatrix} = 2$$ 
and rank $H_N^{(0)} = 2, N \ge 2$. The 2 eigenvalues that we can already
compute, are $\lambda_5^{(0)} = \exp(2\pi\im 55/100)$ and $\lambda_7^{(0)}
= \exp(2\pi\im 45/100)$, satisfying
$$\begin{pmatrix} f_5 & f_{10} \\ f_{10} & f_{15} \end{pmatrix} v = \lambda 
\begin{pmatrix} f_0 & f_5 \\ f_5 & f_{10} \end{pmatrix} v.$$
From the Vandermonde system
$$\begin{pmatrix} 1 & 1 \\ \lambda_5^{(0)} & \lambda_7^{(0)}
\end{pmatrix}
\begin{pmatrix} \alpha_5^{(0)} \\ \alpha_7^{(0)} \end{pmatrix} =
\begin{pmatrix} f_0 \\
f_5 \end{pmatrix}$$
we find $\alpha_5^{(0)} = e^{\im 2\pi 72/100} - e^{\im 2\pi 32/100}$ and
$\alpha_7^{(0)} = 1$.
We now need to ask ourselves
whether $n_0$ truly equals 2 or whether some cancellation of terms has
happened.
With $\rho=12$ we find that
\begin{multline}
f_{5j+12k} = \\ \left( e^{\im 2\pi 12k/100} - e^{\im 2\pi 52k/100} + e^{\im
2\pi 92k/100} - e^{\i 2\pi 32k/100} \right) \exp(2\pi\im 5j/100) + \\
\left( e^{\im 2\pi 72/100}e^{\im 2\pi 32k/100} - e^{\im 2\pi 32/100}e^{\im
2\pi 72k/100} \right) \exp(2\pi\im 55j/100) + \\ 
e^{\im 2\pi 8k/100} \exp(2\pi\im 45j/100).
\end{multline}
For $k=1$ we hit an accidental zero for the coefficient of $\exp(2\pi\im
55j/100)$ and therefore
$$\text{rank } H^{(12)}_3 = \text{rank} \begin{pmatrix} f_{12} & f_{17} &
f_{22} \\ f_{17} & f_{22} & f_{27} \\ f_{22} & f_{27} & f_{32}
\end{pmatrix} = 2$$
again, with rank $H_N^{(12)}=2, N \ge 2$.
The generalized eigenvalues satisfying
$$\begin{pmatrix} f_{17} & f_{22} \\ f_{22} & f_{27} \end{pmatrix} v = \lambda
\begin{pmatrix} f_{12} & f_{17} \\ f_{17} & f_{22} \end{pmatrix} v,$$
namely $\lambda_1^{(0)}=\exp(2\pi\im 5/100)$ and
$\lambda_7^{(0)}=\exp(2\pi\im 45/100)$,
also belong to the $n_0$ eigenvalues that are identifiable from the
evaluations at the multiples of $r\Delta$, since a shift does not change
the generalized eigenvalues, only their coefficients.
From the Vandermonde system
$$\begin{pmatrix} 1 & 1 \\ \lambda_1^{(0)} & \lambda_7^{(0)}
\end{pmatrix}
\begin{pmatrix} \alpha_1^{(1)}(1) \\ \alpha_7^{(1)}(1) \end{pmatrix}
= \begin{pmatrix}
f_{12} \\ f_{17} \end{pmatrix}$$
we find $\alpha_1^{(1)}(1)$ and $\alpha_7^{(1)}(1)$. Apparently $n_0$
equals at least 3, because in the first bunch computed from the samples
$f_{5j+12k}$ with $k=0$ we find two eigenvalues and in the second bunch
with $k=1$ we find one more. Bringing these results together results in
the intermediate estimates 
\begin{align*}
h_1=1&: \alpha_{h_1}^{(1)}(0) = 0, \alpha_{h_1}^{(1)}(1) = e^{\im 2\pi
12/100} - e^{\im 2\pi 52/100} + e^{\im 2\pi 92/100} - e^{\i 2\pi 32/100}, \\
h_2=5&: \alpha_{h_2}^{(1)}(0) = e^{\im 2\pi 72/100} - e^{\im 2\pi 32/100},
\alpha_{h_2}^{(1)}(1) = 0, \\
h_3=7&: \alpha_{h_3}^{(1)}(0) = 1, \alpha_{h_3}^{(1)}(1) = e^{\im 2\pi
8/100}, 
\end{align*}
where the values $h_i$ are merely mentioned as
a guideline and are not explicitly computed.
Remember that in real-life experiments 
the indices $h_i$ are not known and need not be
known. They are revealed as the algorithm progresses.

Let us turn our attention to larger values of $k$ to have the current estimate
$n_0=3$ confirmed and to extract all $n$ distinct terms.
As described in
Section~\ref{sect:four} on the disentangling of collisions, we 
continue sampling at
multiples of the shift, namely we collect the $f_{jr+k\rho}$ for $k>1$.
With $k=2$ we obtain
\begin{multline}
f_{5j+24} = \left( e^{\im 2\pi 24/100} - e^{\im 2\pi 4/100} + e^{\im
2\pi 84/100} - e^{\i 2\pi 64/100} \right) \exp(2\pi\im 5j/100) + \\
\left( e^{\im 2\pi 36/100} - e^{\im 2\pi 76/100} \right) \exp(2\pi\im 55j/100) 
+ \\ e^{\im 2\pi 16/100} \exp(2\pi\im 45j/100)
\end{multline}
and
$$\text{rank } H^{(24)}_4 = \text{rank} \begin{pmatrix} f_{24} & f_{29} &
f_{34} & f_{39} \\ \vdots & \udots & & \vdots \\ f_{39} & \ldots & & f_{54}
\end{pmatrix} = 3$$
with rank $H_N^{(24)}=3, N \ge 4$.
Merely for completeness we compute the generalized eigenvalues satisfying
$$\begin{pmatrix} f_{29} & \ldots & f_{39} \\ \vdots & \udots & \vdots \\ 
f_{39} & \ldots & f_{49} \end{pmatrix} v = \lambda
\begin{pmatrix} f_{24} & \ldots & f_{34} \\ \vdots & \udots & \vdots \\ 
f_{34} & \ldots & f_{44} \end{pmatrix} v.$$
We find $\lambda_1^{(0)}=\exp(2\pi\im 5/100),
\lambda_5^{(0)}=\exp(2\pi\im 55/100), \lambda_7^{(0)}=\exp(2\pi\im
45/100)$, which confirms our earlier obtained combined result. Hence
$n_0=3$.
We also compute  the values for $\alpha_i^{(1)}(2), i=1, 5, 7$
from the Vandermonde system
$$\begin{pmatrix} 1 & 1 & 1 \\ 
\lambda_1^{(0)} & \lambda_5^{(0)} & \lambda_7^{(0)} \\
(\lambda_1^{(0)})^2 & (\lambda_5^{(0)})^2 & (\lambda_7^{(0)})^2 
\end{pmatrix} \begin{pmatrix} \alpha_1^{(1)}(2) \\ \alpha_5^{(1)}(2) \\
\alpha_7^{(1)}(2) \end{pmatrix} = \begin{pmatrix} f_{24} \\ f_{29} \\ f_{34}
\end{pmatrix}.$$
The purpose now is to find
out how many terms are in the expressions $\alpha_i^{(1)}(k)$
for each $i$ retrieved so far. We compute $\alpha_1^{(1)}(k),
\alpha_5^{(1)}(k), \alpha_7^{(1)}(k)$ for $k \ge 3$ from 
$$\begin{pmatrix} 1 & 1 & 1 \\ 
\lambda_1^{(0)} & \lambda_5^{(0)} & \lambda_7^{(0)} \\
(\lambda_1^{(0)})^2 & (\lambda_5^{(0)})^2 & (\lambda_7^{(0)})^2 
\end{pmatrix} \begin{pmatrix} \alpha_1^{(1)}(k) \\ \alpha_5^{(1)}(k) \\
\alpha_7^{(1)}(k) \end{pmatrix} = \begin{pmatrix} f_{12k} \\ f_{5+12k} \\ f_{10+12k}
\end{pmatrix}, \qquad k=3, 4, \ldots,$$
which reuses the $n_0\times n_0$ Vandermonde coefficient matrix from above.

Let us write $k=2\kappa-2$, so that when we increase $\kappa$ by 1 then
$k$ is increased by 2. We check the rank of the $\kappa\times\kappa$
matrices 
$$\begin{pmatrix} \alpha_i^{(1)}(0) & \ldots & \alpha_i^{(1)}(\kappa-1) \\
\vdots & \udots & \vdots \\
\alpha_i^{(1)}(\kappa-1) & \ldots & \alpha_i^{(1)}(2\kappa-2) 
\end{pmatrix}, \qquad i=1,5, 7.$$
When pursuing the shifts up to $k=9$, meaning $\kappa=5$, we find that
for $i=1$ the rank is 4, for $i=5$ the rank is 2 and for $i=7$ the rank is
1, leading to a grand total of $n=7$ distinct terms. We now separate the
terms that are hiding in each of the collisions by computing the
generalized eigenvalues satisfying
\begin{multline}
\begin{pmatrix} \alpha_i^{(1)}(1) & \ldots & \alpha_i^{(1)}(\kappa) \\
\vdots & \udots & \vdots \\
\alpha_i^{(1)}(\kappa) & \ldots & \alpha_i^{(1)}(2\kappa-1) \end{pmatrix} 
v = \lambda
\begin{pmatrix} \alpha_i^{(1)}(0) & \ldots & \alpha_i^{(1)}(\kappa-1) \\
\vdots & \udots & \vdots \\
\alpha_i^{(1)}(\kappa-1) & \ldots & \alpha_i^{(1)}(2\kappa-2)
\end{pmatrix} v, \\ \qquad i=1, 5, 7.
\end{multline}
We find
\begin{align*}
i=1, \kappa=4 : &\lambda_1^{(1)} = \exp(2\pi\im 12/100), \lambda_2^{(1)} =
\exp(2\pi\im 52/100), \\
&\lambda_3^{(1)} = \exp(2\pi\im 92/100),
\lambda_4^{(1)} = \exp(2\pi\im 32/100), \\
i=5, \kappa=2 : &\lambda_5^{(1)} = \exp(2\pi\im 32/100), \lambda_6^{(1)} =
\exp(2\pi\im 72/100), \\
i=7, \kappa=1 : &\lambda_7^{(1)} = \exp(2\pi\im 8/100).
\end{align*}
At this stage we have all the information to reconstruct the non-aliased
generalized eigenvalues:
\begin{align*}
\lambda_1^{(0)}, \lambda_1^{(1)} &\to \lambda_1=\exp(2\pi\im 1/100), \\
\lambda_1^{(0)}, \lambda_2^{(1)} &\to \lambda_2=\exp(2\pi\im 21/100), \\
\lambda_1^{(0)}, \lambda_3^{(1)} &\to \lambda_3=\exp(2\pi\im 41/100), \\
\lambda_1^{(0)}, \lambda_4^{(1)} &\to \lambda_4=\exp(2\pi\im 61/100), \\
\lambda_5^{(0)}, \lambda_5^{(1)} &\to \lambda_5=\exp(2\pi\im 11/100), \\
\lambda_5^{(0)}, \lambda_6^{(1)} &\to \lambda_6=\exp(2\pi\im 31/100), \\
\lambda_7^{(0)}, \lambda_7^{(1)} &\to \lambda_7=\exp(2\pi\im 9/100). 
\end{align*}
Remains to compute the individual linear coefficients of each of the 7
terms. We compute $\alpha_1, \alpha_2, \alpha_3, \alpha_4$
from
$$\begin{pmatrix} 1 & 1 & 1 & 1 \\
\lambda_1^{(1)} & \lambda_2^{(1)} & \lambda_3^{(1)} & \lambda_4^{(1)} \\
(\lambda_1^{(1)})^2 & (\lambda_2^{(1)})^2 & (\lambda_3^{(1)})^2 &
(\lambda_4^{(1)})^2 \\
(\lambda_1^{(1)})^3 & (\lambda_2^{(1)})^3 & (\lambda_3^{(1)})^3 &
(\lambda_4^{(1)})^3 \end{pmatrix} \begin{pmatrix} \alpha_1 \\ \alpha_2 \\ 
\alpha_3 \\ \alpha_4 \end{pmatrix} \begin{pmatrix} \alpha_1^{(1)}(0) \\
\alpha_1^{(1)}(1) \\ \alpha_1^{(1)}(2) \\ \alpha_1^{(1)}(3)
\end{pmatrix},$$
the coefficients $\alpha_5$ and $\alpha_6$ from
$$\begin{pmatrix} 1 & 1 \\
\lambda_5^{(1)} & \lambda_6^{(1)} \end{pmatrix} \begin{pmatrix} 
\alpha_5 \\ \alpha_6 \end{pmatrix} \begin{pmatrix} \alpha_5^{(1)}(0) \\
\alpha_5^{(1)}(1) \end{pmatrix}.$$
The coefficient $\alpha_7$ is given by $\alpha_7=\alpha_7^{(0)}
=\alpha_7^{(1)}(0)$ because there were no collisions in that term.

\subsection{Full algorithm in pseudocode}\label{sect:five.four}

An algorithm covering the eventuality of
the above scenarios reads as follows. We assume that $r>1$ otherwise a
classical Prony analysis applies. 

So far we used the notation $n$ for the number of exponential terms in the
signal, which we often don't know up front. Moreover, the data are usually
noisy, so that it is best to add another number of terms in order to model the
noise. We denoted the latter in the previous sections
by $N-n$ so that the total number of terms we want to identify accumulates 
to $N$. To this end at least $2N$ samples are required, even without
breaking the Shannon-Nyquist rate. We denote the
number of samples collected at the uniformly distributed points $t_{jr}$
by the number $M \ge 2N$. These allow us to build the square Hankel matrices
$H^{(0)}_N$ and $H^{(1)}_N$ or somewhat larger rectangular $(M-N)\times N$
versions of these matrices. When sampling at the shifted locations
$t_{jr+k\rho}$ we collect for each $k$ not $M$ but $m$ samples 
where $\lfloor m/2
\rfloor > n$. Using the latter we can build the Hankel matrices
$H^{(\rho)}_{\lfloor m/2 \rfloor}$. 
Often the total number of samples and the amount of
undersampling are dictated by the circumstances and the constraints under
which the analysis is performed. 

We emphasize that $n_0$ indicates the number of terms in the exponential
sum after possible collisions, including the vanished ones due to
cancellation in the coefficients. Also, the time step
$\Delta \in \rz$ satisfies $\Delta \le 1/\Omega$.
With this in mind the algorithm continues as follows.
\medskip

{\bf Algorithm.}
\medskip

{\sl Input bounds on $n$, subsampling factor $r$ and shift term $\rho$:}

\begin{itemize}
\item $M, N, m \in \nz$ with $M \ge 2N, N \ge n, \lfloor m/2 \rfloor > n$.
\item $r \in \nz, \rho \in \gz$ with $r > 1$ and $\gcd(r, \rho)=1$.
\end{itemize}
\medskip

{\sl A0. Obtain $n_0, \lambda_i^{(0)}, \alpha_i^{(0)}$:}

\begin{itemize}
\item Collect the samples $f_{jr} = \phi(t_{jr}), j=0, \ldots, M-1$ and
estimate $n_0\le n$ by the numerical rank of the matrix $H_N^{(0)}$. 

\item For one or more $1 \le k \le 2n-1$ collect the samples $f_{jr+k\rho}= 
\phi(t_{jr+k\rho}), j=0, \ldots, m-1$
and compute the numerical rank $n_k$ of $H_{\lfloor m/2 \rfloor}^{(k\rho)}$. 

\item From these different views on the number of collided terms in
the exponential sum, we find that the correct value for $n_0$ is
$n_0=\max_k n_k$.
  
\item Compute for $i=1, \ldots, n_0$ the generalized eigenvalues 
$\lambda_i^{(0)}$ and the coefficients $\alpha_i^{(0)}$ as in Example 5.3.
  
\item Either $N\times N$ Hankel and $2N \times N$
Vandermonde systems are used or their
least squares $(M-N)\times N$ and $M\times N$ versions.
\end{itemize}
\medskip  

{\sl A1. Obtain $\alpha_i^{(1)}(k)$ and $h_{i+1}-h_i$
for $i=1, \ldots, n_0$ and $1 \le k \le K < 2n$:}

Put $\alpha_i^{(1)}(0):=\alpha_i^{(0)}, k=1$ and execute the for loop:
\begin{enumerate}

\item compute the $\alpha_i^{(1)}(k)$ from \eqref{step1} or its $m\times
n_0$ least squares version,
  
\item collect or reuse the samples $f_{jr+(k+1)\rho} = \phi(t_{jr+(k+1)\rho}), j=0, \ldots,
m-1$,

\item compute the $\alpha_i^{(1)}(k+1)$ from \eqref{step1} or its $m\times
n_0$ least squares version,

\item compute the numerical rank $\nu_i(\kappa)$ of the
$(\kappa+1)\times(\kappa+1)$ matrix 
$$\begin{bmatrix}
\alpha_i^{(1)}(0) & \cdots & \alpha_i^{(1)}(\kappa) \\
\vdots & & \vdots \\
\alpha_i^{(1)}(\kappa) & \cdots & \alpha_i^{(1)}(k+1) 
\end{bmatrix}, \qquad 2\kappa = k+1,$$
   
\item if $\nu_i(\kappa)= \nu_i(\kappa-1)$: 
\begin{itemize}
\item then $h_{i+1}-h_i = \kappa$,
\item else $k:= k+2$, collect or reuse the samples $f_{jr+k\rho}, j=0, \ldots,
m-1$ and goto 1.
\end{itemize}

\item compute the generalized eigenvalues $\lambda_\ell^{(1)}, \ell=h_i,
\ldots, h_{i+1}-1$ in \eqref{four.three}, 

\item compute the $\alpha_\ell, \ell=h_i, \ldots, h_{i+1}-1$ from 
\eqref{alpha_i1}.
\end{enumerate}

End the for loop.
\medskip

{\sl Output number of terms $n$ and the parameters $\phi_i, \alpha_i$:}

From
\begin{itemize}
\item $\lambda_\ell^{(1)}, \ell=h_1, \ldots, 
h_{n_0+1}-1$ with $h_1=1, h_{n_0+1}-1=n$ 
\item and $\lambda_\ell^{(0)}$ with
$\lambda_{h_i}^{(0)} = \cdots =
\lambda_{h_{i+1}-1}^{(0)}, i=1, \ldots, n_0$ 
\end{itemize}
the $\phi_i, i=1, \ldots, n$ can be recovered.

The $\alpha_i, i=1, \ldots, n$ are computed from \eqref{alpha_i1} as in
\eqref{qedemo}.

\begin{table}[h]
\centering
\begin{tabular}{|c|c|c|}
\hline
$i$ & $\alpha_i$ & $\phi_i$ \\ \hline 
1 & $6.5\exp(0.15\im)$ &$-0.19-\im 2 \pi 453.1$ 
\\
2 & $6.8$ &$-0.132-\im 2 \pi 452.19$   
\\
3 & $6.8\exp(0.3\im)$ &$-0.183-\im 2 \pi 451.02$
\\
4 & $6.4\exp(0.9\im)$ &$-0.11-\im 2 \pi 450.21$ 
\\
5 & $7.1\exp(0.7\im)$ & $-0.21-\im 2 \pi 448.39$ 
\\ \hline 
6 & $4.71\exp(0.12\im)$ & $-0.106 -\im 2 \pi 132.5$
\\
7 & $3.9\exp(0.1\im)$ & $-0.129 -\im 2 \pi 131.4$
\\
8 & $7.2\exp(-0.234\im)$ & $-0.203 -\im 2 \pi
130.01$  
\\
9 & $7.43\exp(0.2\im)$ & $-0.16 -\im 2 \pi 129.17$
\\
10 & $4.4\exp(-0.52\im)$ & $-0.19 -\im 2 \pi
128.39$
\\ \hline 
11 & $3\exp(0.21\im)$ & $-0.101+ \im 2 \pi 9.1$
\\
12 & $3\exp(-0.8\im)$ & $-0.127 + \im 2 \pi
11.81$ 
\\ \hline
13 & $7.2\exp(-0.106\im)$ & $-0.21 + \im 2 \pi
126.01$ 
\\
14 & $6.53\exp(0.2\im)$ & $-0.15 + \im 2 \pi
127.62$
\\
15 & $6.7\exp(-0.3\im)$ & $-0.173 + \im 2 \pi
128.98$
\\ \hline 
16 & $6.8\exp(-0.15\im)$ & $-0.11 + \im 2 \pi
334.01$ 
\\
17 & $6\exp(0.26\im)$ & $-0.12  + \im 2 \pi
335.18$ 
\\
18 & $7.1\exp(-0.2\im)$ & $-0.157 + \im 2 \pi
336.01$ 
\\
19 & $7.1$ & $-0.120 + \im 2 \pi 337.91$ 
\\
20 & $6\exp(-0.1\im)$ &$-0.18 + \im 2 \pi
339.61$ \\ \hline
\end{tabular}
\caption{Collision-free example.}\label{table:one}
\end{table}
\begin{table}[h]
\centering
\begin{tabular}{|c|c|c|}
\hline
$i$ & $\alpha_i$ & $\phi_i$ \\ \hline 
1 & $18$ & $\im 2 \pi 191.9$ \\
2 & $-20$ & $\im 2 \pi 291.9$ \\
3 & $20$ & $\im 2 \pi 391.9$ \\ \hline 
4 & $5$ & $\im 2 \pi 526.2$ \\ \hline 
5 & $5$ & $\im 2 \pi 858.1$ \\
6 & $11$ & $\im 2 \pi 958.1$  \\ \hline
\end{tabular}
\caption{Example where collisions occur.}\label{table:two}
\end{table}


\section*{Acknowledgements}

The authors sincerely thank Engelbert Tijskens of the Universiteit Antwerpen 
for making the documented Matlab code available that is downloadable from the
webpage {\tt cma.uantwerpen.be/publications} and that allows the reader to
rerun all the numerical illustrations included in the paper and even
several variations thereof. 

This work was partially supported by a Research Grant of the FWO-Flanders 
(Flemish Science Foundation) and a Proof of Concept project of the 
University of Antwerp (Belgium).



\begin{thebibliography}{10}
\expandafter\ifx\csname url\endcsname\relax
  \def\url#1{\texttt{#1}}\fi
\expandafter\ifx\csname urlprefix\endcsname\relax\def\urlprefix{URL }\fi
\expandafter\ifx\csname href\endcsname\relax
  \def\href#1#2{#2} \def\path#1{#1}\fi

\bibitem{CandesFernandez:2014}
E.~J.~Cand\`{e}s, C.~Fernandez-Granda, Towards a mathematical theory of
  super-resolution, Communications on Pure and Applied Mathematics 67~(6)
  (2014) 906--956.

\bibitem{Moitra:2015}
A.~Moitra, Super-resolution, extremal functions and the condition number of
  Vandermonde matrices, in: Proceedings of the Forty-seventh Annual ACM
  Symposium on Theory of Computing, STOC '15, ACM, 2015, pp. 821--830.

\bibitem{Kammler}
D.~W. Kammler, Approximation with sums of exponentials in $l_p[0, \infty)$,
  Journal of Approximation Theory 16~(4) (1976) 384--408.

\bibitem{Va:fit:85}
J.~Varah, On fitting exponentials by nonlinear least squares, SIAM Journal on 
Scientific and Statistical Computing 6~(1) (1985) 30--44.

\bibitem{Ha.Ka:eff:97}
B.~Halder, T.~Kailath, Efficient estimation of closely spaced sinusoidal
  frequencies using subspace-based methods, {IEEE} Signal Processing Letters
  4~(2) (1997) 49--51.

\bibitem{Sc:mul:86}
R.~Schmidt, Multiple emitter location and signal parameter estimation, IEEE
  Transactions on Antennas and Propagation 34~(3) (1986) 276--280.

\bibitem{Ro.Ka:esp:89}
R.~Roy, T.~Kailath, {ESPRIT}-estimation of signal parameters via rotational
  invariance techniques, {IEEE} Transactions on Acoustics, Speech, and Signal 
  Processing 37~(7) (1989) 984--995.

\bibitem{Hu.Sa:mat:90}
Y.~Hua, T.~K. Sarkar, Matrix pencil method for estimating parameters of
  exponentially damped/undamped sinusoids in noise, IEEE Transactions on 
  Acoustics, Speech, and Signal Processing 38 (1990) 814--824.

\bibitem{Go.Mi.ea:sta:99}
G.~Golub, P.~Milanfar, J.~Varah, A stable numerical method for inverting shape
  from moments, SIAM Journal on Scientific Computing 21 (1999) 1222--1243.

\bibitem{Neumaier}
S.~Das, A.~Neumaier, Solving overdetermined eigenvalue problems, {SIAM} Journal
  on Scientific Computing 35~(2) (2013) A541--A560.

\bibitem{Be.Mo:app:05}
G.~Beylkin, L.~Monz{\'o}n, On approximation of functions by exponential sums,
  Applied and Computational Harmonic Analysis 19~(1) (2005) 17--48.

\bibitem{Po.Ta:par:10}
D.~Potts, M.~Tasche, Parameter estimation for exponential sums by approximate
  {P}rony method, Signal Processing 90 (2010) 1631--1642.

\bibitem{Po.Ta:par:13}
D.~Potts, M.~Tasche, Parameter estimation for nonincreasing exponential sums by
  {P}rony-like methods, Linear Algebra and its Applications 439~(4) (2013) 1024--1039.

\bibitem{Ny:cer:28}
H.~Nyquist, Certain topics in telegraph transmission theory,  Transactions of the 
American Institute of Electrical Engineers 47~(2) (1928) 617--644.

\bibitem{Sh:com:49}
C.~E. Shannon, Communication in the presence of noise,  Proceedings of the IRE 
37 (1949) 10--21.

\bibitem{Ca.Ro.ea:rob:06}
E.~J. Cand{\`e}s, J.~Romberg, T.~Tao, Robust uncertainty principles: exact
  signal reconstruction from highly incomplete frequency information, {IEEE}
  Transactions on Information Theory 52~(2) (2006) 489--509.

\bibitem{Donoho:2006}
D.~L. Donoho, Compressed sensing, {IEEE} Transactions on Information Theory
  52~(4) (2006) 1289--1306.

\bibitem{VMB02}
M.~Vetterli, P.~Marziliano, T.~Blu, Sampling signals with finite rate of
  innovation, {IEEE} Transactions on Signal Processing 50~(6) (2002)
  1417--1428.

\bibitem{VaPa2011-2}
P.~P. Vaidyanathan, P.~Pal, Sparse sensing with co-prime samplers and arrays,
  {IEEE} Transactions on Signal Processing 59~(2) (2011) 573--586.

\bibitem{TEN14}
Z.~Tan, Y.~C. Eldar, A.~Nehorai, Direction of arrival estimation using co-prime
  arrays: A super resolution viewpoint, {IEEE} Transactions on Signal
  Processing 62~(21) (2014) 5565--5576.

\bibitem{Cu.Le:sma:12*b}
A.~Cuyt, W.-s. Lee, Smart data sampling and data reconstruction, patent
  EP2745404B1. 

\bibitem{Cu.Le:sma:12}
A.~Cuyt, W.-s. Lee, Smart data sampling and data reconstruction, patent US
  9,690,740. 

\bibitem{Pl.Ta:pro:14}
G.~Plonka, M.~Tasche, {P}rony methods for recovery of structured functions,
  GAMM-Mitt. 37~(2) (2014) 239--258.

\bibitem{He:app:74}
P.~Henrici, Applied and computational complex analysis {I}, John Wiley \& Sons,
  New York, 1974.

\bibitem{Ka.Le:ear:03}
E.~Kaltofen, W.-s. Lee, Early termination in sparse interpolation algorithms,
  Journal of Symbolic Computation 36~(3-4) (2003) 365--400. 

\bibitem{Cu.Ts.ea:fai:18}
A.~Cuyt, M.~Tsai, M.~Verhoye, W.-s. Lee, Faint and clustered components in
  exponential analysis, Applied Mathematics and Computation 327 (2018) 93--103.

\bibitem{Be.Go.ea:num:07}
B.~Beckermann, G.~Golub, G.~Labahn, On the numerical condition of a generalized
  {H}ankel eigenvalue problem, Numerische Mathematik 106~(1) (2007) 41--68.

\bibitem{Ga:nor:75}
W.~Gautschi, Norm estimates for inverses of {V}andermonde matrices, Numerische 
  Mathematik 23 (1975) 337--347.

\bibitem{Cu.Le:ana:17}
A.~Cuyt, W.-s. Lee, An analog {C}hinese {R}emainder {T}heorem, Tech. rep.,
  Universiteit Antwerpen (2017).

\end{thebibliography}
\end{document}